\documentclass[11pt,leqno]{article}
\normalfont
\renewcommand{\baselinestretch}{1.1735}
\usepackage{amsfonts}

\newfont{\eulercursive}{eurm10 at 11pt}
\newcommand{\myl}{\mbox{\eulercursive `}}

\newcommand{\velt}{\mathbf{v}}
\newcommand{\uelt}{\mathbf{u}}

\newcommand{\QED}{\raisebox{0.5mm}{\fbox{\rule{0mm}{1.5mm}\ }}}

\newcommand{\BBProposition}{Proposition 2.1} 
\newcommand{\TwoGeneratorAnalysis}{Lemma 3.1}
\newcommand{\HumphreysTheorem}{Theorem 3.2}
\newcommand{\HumphreysCorollary}{Corollary 3.3} 
\newcommand{\OneToOneLemma}{Lemma 3.4}
\newcommand{\ReducedLemma}{Lemma 3.5}
\newcommand{\TFAE}{Theorem 3.6}
\newcommand{\TFAECorollary}{Corollary 3.7}
\newcommand{\TFAECorollaryTheorems}{Theorems 3.2 and 3.6} 
\newcommand{\PositiveToNegativeRootsLemma}{Lemma 3.8}
\newcommand{\LengthProposition}{Theorem 3.9}
\newcommand{\ListForCorollary}{Theorems 3.6 and 3.9}
\newcommand{\UnitalCorollary}{Corollary 3.10}
\newcommand{\LengthCorollary}{Corollary 3.11}
\newcommand{\UnitalExample}{Example 3.12}
\newcommand{\BrinkHowlettTheorem}{Theorem 3.13}
\newcommand{\UnitalExampleFigure}{Figure 3.1}
\newcommand{\TitsConeConvergenceResult}{Proposition 4.1}
\newcommand{\DeodharProp}{Lemma 4.2}
\newcommand{\HRTProofResults}{Corollary 3.10 and Lemma 4.2}
\newcommand{\HRTResult}{Proposition 4.3}
\newcommand{\TheoremList}{Propositions 4.1 and 4.3}
\newcommand{\TitsConeProp}{Proposition 4.4}

\newcommand{\TitsConeFiniteResult}{Theorem 4.5}
\newcommand{\KumarRemark}{Remark 4.6}
\newcommand{\VinbergRemark}{Remark 4.7}
\newcommand{\TitsConeExample}{Example 4.8}

\newcommand{\disjointunion}{\setlength{\unitlength}{0.14cm}
\ 
\begin{picture}(2,2) 
\put(0,0){$\cup$}
\put(0.9675,1.5){\circle*{0.5}}
\end{picture}\ }

\newcommand{\CircleInteger}[1]{
\setlength{\unitlength}{0.14cm}
\begin{picture}(3,2) 
\put(1,1){\circle{2}}
\put(0.55,0.6){{\tiny #1}}
\end{picture}
}

\newcommand{\CircleIntegerm}{
\setlength{\unitlength}{0.14cm}
\begin{picture}(3,2) 
\put(1,1){\circle{2}}
\put(0.3,0.6){{\tiny $m$}}
\end{picture}
}

\newcommand{\CircleInfty}{
\setlength{\unitlength}{0.14cm}
\begin{picture}(3,2) 
\put(1,1){\circle{2}}
\put(0.1,0.6){{\tiny $\infty$}}
\end{picture}
}

\newcommand{\TwoCitiesGraphWithLabels}{
\setlength{\unitlength}{0.75in}
\begin{picture}(1.65,0.25)
\put(0.25,0){\begin{picture}(1,0)
            \put(0,0.1){\circle*{0.05}}
            \put(-0.20,-0.05){\large $\gamma_{i}$}
            \put(1,0.1){\circle*{0.05}}
            \put(1.05,-0.05){\large $\gamma_{j}$}
            \put(0,0.1){\line(1,0){1}}
            \put(0.2,0.1){\vector(1,0){0.1}}
            \put(0.8,0.1){\vector(-1,0){0.1}}
            \put(0.225,-0.05){\footnotesize $p$}
            \put(0.71,-0.05){\footnotesize $q$}
            \end{picture}}
\end{picture}}

\newcommand{\TwoCitiesGraphWithoutLabels}{
\setlength{\unitlength}{0.75in}
\begin{picture}(1.65,0.25)
\put(0.25,0){\begin{picture}(1,0)
            \put(0,0.1){\circle*{0.05}}
            \put(1,0.1){\circle*{0.05}}
            \put(0,0.1){\line(1,0){1}}
            \put(0.2,0.1){\vector(1,0){0.1}}
            \put(0.8,0.1){\vector(-1,0){0.1}}
            \put(0.225,-0.05){\footnotesize $p$}
            \put(0.71,-0.05){\footnotesize $q$}
            \end{picture}}
\end{picture}}

\newcommand{\BowTie}{\setlength{\unitlength}{1.55in}
\begin{picture}(2.5,1.2)
\put(0,0){
\begin{picture}(0,1.2)
            \put(0.8,0.6){\line(1,0){0.81}}
            \put(1.1,0.625){\CircleInfty}
            \put(0.8,0.6){\vector(1,0){0.2}}
            \put(1.6,0.6){\vector(-1,0){0.2}}
            \put(0.95,0.51){\footnotesize $2$}
            \put(1.4,0.51){\footnotesize $2$}
            \put(0.3,0.1){\line(0,1){1}}
            \put(0.3,0.1){\line(1,1){0.5}}
            \put(0.3,1.1){\line(1,-1){0.5}}
            \put(0.8,0.6){\circle*{0.05}}
            \put(0.3,1.1){\circle*{0.05}}
            \put(0.3,0.1){\circle*{0.05}}
            \put(0.85,0.675){\small $\gamma_{2}$}
            \put(0.125,0.05){\small $\gamma_{3}$}
            \put(0.125,1.1){\small $\gamma_{1}$}
            \put(0.3,0.3){\vector(0,1){0.1}}
            \put(0.195,0.35){\footnotesize $1$}
            \put(0.3,0.9){\vector(0,-1){0.1}}
            \put(0.195,0.8){\footnotesize $3$}
            \put(0.255,0.575){\CircleInteger{6}}
            \put(0.3,0.1){\vector(1,1){0.15}}
            \put(0.8,0.6){\vector(-1,-1){0.15}}
            \put(0.45,0.15){\footnotesize $2$}
            \put(0.675,0.375){\footnotesize $1$}
            \put(0.425,0.375){\CircleInteger{4}}
            \put(0.3,1.1){\vector(1,-1){0.15}}
            \put(0.8,0.6){\vector(-1,1){0.15}}
            \put(0.425,0.755){\CircleInteger{3}}
            \put(0.435,1){\footnotesize $1/5$}
            \put(0.69,0.75){\footnotesize $5$}
\end{picture}
}
\put(1.6,0){
\begin{picture}(0.6,1.2)
            \put(0,0.6){\circle*{0.05}}
            \put(0.5,0.1){\circle*{0.05}}
            \put(0.5,1.1){\circle*{0.05}}
            \put(0,0.6){\line(1,1){0.5}}
            \put(0,0.6){\line(1,-1){0.5}}
            \put(0.5,0.1){\line(0,1){1}}
            \put(0.5,0.3){\vector(0,1){0.1}}
            \put(0.55,0.35){\footnotesize $1+\sqrt{5}$}
            \put(0.5,0.9){\vector(0,-1){0.1}}
            \put(0.55,0.8){\footnotesize $(1+\sqrt{5})/4$}
            \put(0,0.6){\vector(1,1){0.15}}
            \put(0.5,1.1){\vector(-1,-1){0.15}}
            \put(0.075,0.75){\footnotesize $7$}
            \put(0.25,1){\footnotesize $1/7$}
            \put(0,0.6){\vector(1,-1){0.15}}
            \put(0.5,0.1){\vector(-1,1){0.15}}
            \put(0,0.375){\footnotesize $7/2$}
            \put(0.25,0.15){\footnotesize $2/7$} 
            \put(-0.125,0.675){\small $\gamma_{5}$}
            \put(0.58,0.05){\small $\gamma_{6}$}
            \put(0.58,1.1){\small $\gamma_{4}$}
            \put(0.345,0.575){\CircleInteger{5}}
            \put(0.175,0.375){\CircleInteger{3}}
            \put(0.175,0.755){\CircleInteger{3}}
\end{picture}
}
\end{picture}
}

\begin{document}

\newpage
\setcounter{page}{1} 
\renewcommand{\baselinestretch}{1}

\vspace*{-0.7in}
\hfill {\footnotesize December 30, 2009}

\begin{center}
{\large \bf Root systems for asymmetric geometric representations of 
Coxeter groups}

\renewcommand{\thefootnote}{1}
Robert G.\ Donnelly\footnote{Email: 
{\tt rob.donnelly@murraystate.edu},\ Fax: 
1-270-809-2314}

\vspace*{-0.075in}
Department of Mathematics and Statistics, Murray State
University, Murray, KY 42071
\end{center}

\vspace*{-0.25in}
\begin{abstract}
Results are obtained concerning root systems for asymmetric geometric 
representations of Coxeter groups.  These representations were 
independently introduced by Vinberg and Eriksson, and generalize the 
standard geometric representation of a Coxeter group in such a way as 
to include all Kac--Moody Weyl groups. 
In particular, a characterization of when a non-trivial multiple of a 
root may also be a root is given in the general context.  
Characterizations of when the number of such multiples of a root is 
finite and when the number of positive roots sent to 
negative roots by a group element is finite are also given. 
These characterizations are stated in terms of combinatorial 
conditions on a graph closely related to the Coxeter graph for the 
group. Other finiteness results for the symmetric case  
which are connected to the Tits cone and to a natural 
partial order on positive roots are extended to this asymmetric 
setting. 
\begin{center}


{\small \bf Keywords:}\ {\small Coxeter group, 
geometric representation, root system, Tits cone,\\ Kac--Moody algebra, 
numbers game} 
\end{center}
\end{abstract}

\vspace*{-0.1in}
{\bf \S 1 Introduction.}  
A certain natural symmetric bilinear form is used 
to define the familiar geometric representation of a given 
Coxeter group, often called the ``standard'' geometric representation.  
See 
\cite{Bourbaki} Ch.\ 5, \cite{HumCoxeter} Ch.\ 5, or \cite{BB} \S 4.4.  
These representations are well understood 
and are useful for 
studying Coxeter groups and their applications in many different 
contexts.  
See for example \cite{Gunnells} and references therein.  
Following work of Vinberg and 
Eriksson, when considering geometric representations of Coxeter 
groups in Chapter 4 of the book \cite{BB}, Bj\"{o}rner and Brenti 
initially do not require that the bilinear form be symmetric.  
The purpose here is to 
further study the root systems associated to such representations. 
Much of what we record here 
generalizes the standard theory as presented for example in 
\S 5.3, 5.4, 5.6, and 5.13 of \cite{HumCoxeter} and 
extends \S 4.3 of \cite{BB}. 
Since the form is no longer required to be symmetric, all statements 
here may be applied to the sets of real roots of Kac--Moody algebras.  
This yields new proofs of standard Kac--Moody results (one 
direction of the first statement in  
\TFAECorollary, one direction of the second statement in 
\UnitalCorollary). 

These asymmetric geometric realizations of Coxeter groups were 
introduced by Vinberg in \cite{Vinberg}, for geometric reasons.  
A main focus 
of Vinberg's study is the behavior of the ``fundamental chamber'' (a 
convex polyhedral cone) under the group action. 
In a different context, Lusztig used 
such asymmetric forms when constructing  
certain irreducible representations of Hecke algebras \cite{Lusztig}.  
Eriksson applied asymmetric geometric representations of Coxeter groups 
in \cite{ErikssonThesis} (\S 4.3, \S 6.9, Ch.\ 8) 
and \cite{ErikssonDiscrete} (\S 3, 4) in connection 
with the combinatorial 
numbers game of Mozes \cite{Mozes}. 
While the numbers game is of combinatorial interest in its own right, 
it is also helpful for facilitating computations with Coxeter groups 
and their geometric representations  
(e.g.\ computing orbits, solving the word problem, or finding reduced 
decompositions) and for obtaining combinatorial models of Coxeter 
groups.  See for example \S 4.3 of \cite{BB}.   
The results of this paper 
are needed for our further study of the numbers game in 
\cite{DonNumbers}. 
There we further investigate connections between moves of the game 
and reduced decompositions for group elements, characterize ``full 
commutativity'' of group elements in terms of the game, characterize 
when all positive roots can be obtained from game play, and 
obtain a new Dynkin diagram classification theorem whose 
answer consists of versions of Coxeter graphs for finite Coxeter 
groups. 

The possible asymmetry of the bilinear forms here leads 
to some curious differences with the standard case.   
In Exercise 4.9 of \cite{BB}, the authors point out that without  
symmetry of the bilinear forms, some important properties of root systems 
would not be true. However, we will see that these properties do not 
fail too badly, at least not all of the time. 
In particular, 
we determine precisely when non-trivial scalar multiples of  
roots can also be roots (\HumphreysTheorem), and we 
relate the finiteness of this set of root multiples to 
a combinatorial condition on a graph closely related to the Coxeter graph for 
the group (\TFAE).  
Further, we determine when the number of positive roots sent to 
negative roots by a given group element is finite,  
and we say how this quantity is related to the length of 
the given group element (\LengthProposition).  
An asymmetric version of 
Brink and Howlett's fundamental result on the finiteness of the set of 
``dominance-minimal'' roots is obtained in 
\BrinkHowlettTheorem. 
In \TitsConeFiniteResult, we show that finiteness of an irreducible 
Coxeter group 
is equivalent to certain conditions on the 
asymmetric version of the 
Tits cone. 

The original version of this paper was written with only the numbers 
game motivations above in mind.  Recently, 
for unrelated reasons Proctor decided to relate the treatment of Weyl 
groups in \cite{Kac} and 
\cite{Kumar} to the study of asymmetric geometric representations 
of Coxeter groups in \cite{BB}.  This led to the definition of 
`real Weyl groups' in \cite{ProctorCoxeter} 
and his realization that our \HumphreysTheorem\ 
would play a key role in those notes.  Quoting from an earlier 
draft of \cite{ProctorCoxeter}:  
``There are many statements concerning Weyl groups and the `real' 
roots of Kac--Moody algebras which can at least be conjectured in 
the general context of real Weyl groups.  If still true, it would 
seem that each of these statements should be provable without any 
reference to Lie brackets or to root spaces, if one could formulate 
suitable sufficient conditions for them in terms of real Weyl group 
concepts.  One example of such a statement is ``no `non-trivial' real 
multiple of a real root is also a root''.  Within the general context, 
two successive restricting assumptions (which are both automatically 
satisfied by Weyl groups) guarantee [via our \HumphreysTheorem] that 
this example statement holds true in a context which is still much 
more general than that of Weyl groups or of Section 4.4 of \cite{BB}.''  

At the end of Section 2 we observe that any Kac--Moody Weyl group 
arises as one of our representing groups  $\sigma(W) \subset 
GL(V)$.  
Hence all of our results pertain to the special 
case consisting of arbitrary Kac--Moody Weyl groups. 
Our complete characterizations of the ``no non-trivial 
multiple of a (real) root is also a root'' (\TFAECorollary) and the ``set of 
positive (real) 
roots sent to negative is finite'' (\UnitalCorollary) 
properties are given 
proofs which are naturally set in a general environment which 
encompasses both the standard geometric representations of Coxeter 
groups and Kac--Moody Weyl groups.  Only combinatorial positivity 
arguments are used in these proofs; no references to Lie brackets or 
root spaces are needed. 

{\bf \S 2 Definitions and preliminaries.} 
In this section we present the main objects of interest for 
this paper. 
The crucial information identifying an asymmetric geometric 
representation of a Coxeter group is a certain real matrix analog of a 
generalized Cartan matrix.  We take this matrix as our starting point. 
Fix a positive integer $n$ and a totally ordered set $I_{n}$ with $n$ 
elements (usually $I_{n} := \{1<\ldots<n\}$).  An 
{\em E-generalized Cartan matrix}  
\renewcommand{\thefootnote}{2}
({\em E-GCM})\footnote{Motivation for terminology: 
E-GCM's with integer entries are 
generalizations of `generalized' Cartan matrices (GCM's), which 
are the starting point for the study of 
Kac-Moody algebras. 
Here we use the 
modifier ``E'' because of 
the relationship between these matrices and the combinatorics of 
Eriksson's  
E-games \cite{ErikssonThesis}, 
\cite{ErikssonDiscrete}.}   
 is  
an $n \times n$ matrix $A = (a_{ij})_{i,j \in I_{n}}$  
with real entries satisfying the requirements that each 
main diagonal matrix entry is 2, that all other matrix entries are 
nonpositive, that if 
$a_{ij}$ is nonzero then $a_{ji}$ is also nonzero, and that for $i \not= j$  
either $a_{ij}a_{ji} \geq 4$ or 
$a_{ij}a_{ji} = 4\cos^{2}(\pi/k_{ij})$ for some integer $k_{ij} \geq 
2$. The peculiar quantities $4\cos^{2}(\pi/k)$ appear in 
the developments of \cite{Bourbaki}, \cite{HumCoxeter} 
as the products of transpose entries of a 
symmetric 
matrix for the defining bilinear form of the standard geometric 
representation of a Coxeter group.  
To an $n \times n$ E-generalized Cartan matrix 
$A = (a_{ij})_{i,j \in I_{n}}$ we associate a finite 
graph $\Gamma$ 
as follows:     
The nodes $(\gamma_{i})_{i \in I_{n}}$ of $\Gamma$ are indexed 
by the set $I_{n}$, 
and   an edge is placed between nodes $\gamma_{i}$ and $\gamma_{j}$ 
if and only if $i \not= j$ 
and the matrix entries $a_{ij}$ and $a_{ji}$ are nonzero.  
We display this edge as 
\TwoCitiesGraphWithLabels, 
where $p = -a_{ij}$ and $q = -a_{ji}$.  
We call the pair 
$(\Gamma,A)$ an {\em E-GCM graph}. 
See \UnitalExampleFigure\ for a six-node example. 

Define the associated Coxeter group 
$W(\Gamma,A)$ to be the Coxeter group with identity 
$\varepsilon$, 
generators $\{s_{i}\}_{i \in I_{n}}$, and defining relations $s_{i}^{2} = 
\varepsilon$ for $i \in I_{n}$ and 
$(s_{i}s_{j})^{m_{ij}} = \varepsilon$ for all $i \not= j$, where the 
$m_{ij}$ are determined by: 

\vspace*{-0.1in}
\[m_{ij} = \left\{\begin{array}{cl}
k_{ij} & \mbox{\hspace*{0.25in} if 
$a_{ij}a_{ji} = 4\cos^{2}(\pi/k_{ij})$ for some integer $k_{ij} \geq 2$}\\ 
\infty & \mbox{\hspace*{0.25in} if 
$a_{ij}a_{ji} \geq 4$} 
\end{array}\right.\]  
(Conventionally, $m_{ij} = \infty$ means there is no relation between 
generators 
$s_{i}$ and $s_{j}$.) 
When $A$ is a generalized Cartan matrix or GCM 
(i.e.\ an E-GCM with integer entries), 
then $W(\Gamma,A)$ is a {\em Weyl group}.  In this case, 
$m_{ij}$ is finite only for the pairs 
$\{-a_{ij},-a_{ji}\} = \{0,0\}, \{1,1\}, \{1,2\}, 
\{1,3\}$; the corresponding values of such $m_{ij}$ are $2, 3, 4, 6$. 
One can think of the E-GCM graph as a refinement of the information 
from the Coxeter graph for the associated Coxeter group.   
Observe 
that any Coxeter group on a finite set of generators is isomorphic 
to $W(\Gamma,A)$ for some E-GCM graph $(\Gamma,A)$.  
We let ${\myl}$ denote the length function for $W = W(\Gamma,A)$. 
An expression $s_{i_{1}}s_{i_{2}}{\cdots}s_{i_{p}}$ for an element 
of $W$ is {\em reduced} 
if $\myl(s_{i_{1}}s_{i_{2}}{\cdots}s_{i_{p}}) = p$. 
For $J \subseteq I_{n}$, let 
$W_{J}$ be the subgroup generated by 
$\{s_{i}\}_{i \in J}$, a {\em parabolic} subgroup, and $W^{J} := \{w 
\in W\, |\, \myl(ws_{j}) > \myl(w) \mbox{ for all } j \in J\}$  
is the set of {\em minimal coset representatives}. 
If $J = \{i,j\}$, then $W_{J}$ is a dihedral group of order 
$2m_{ij}$.  

From here on, fix an arbitrary 
E-GCM graph $(\Gamma,A)$ with index set $I_{n}$ 
and associated Coxeter group $W = W(\Gamma,A)$. 
We now define the representations of $W$ 
which are of interest to us here, cf.\ \S 4.2 of \cite{BB}. To fix 
notation that will help set up some subsequent arguments, we present 
some of the details here. 
Let $V$ be a real $n$-dimensional vector space freely generated by 
$(\alpha_{i})_{i \in I_{n}}$. (Elements of this ordered basis are {\em 
simple roots}.)   
Equip $V$ with a possibly asymmetric bilinear form $B: V \times V 
\rightarrow \mathbb{R}$ defined on the 
basis $(\alpha_{i})_{i \in I_{n}}$ by 
$B(\alpha_{i},\alpha_{j}) := \frac{1}{2}a_{ij}$. 
For each $i \in I_{n}$ define an operator 
$S_{i}: V \rightarrow V$ by the rule $S_{i}(v) := v - 
2B(\alpha_{i},v)\alpha_{i}$ for each $v \in V$.  One can check 
that $S_{i}^{2}$ is the identity 
transformation, so $S_{i} \in GL(V)$.  Fix $i \not= j$ and set $V_{i,j} := 
\mathrm{span}_{\mathbb{R}}\{\alpha_{i},\alpha_{j}\}$.  
Observe that $S_{k}(V_{i,j}) \subseteq V_{i,j}$ for $k = i,j$.  Let 
$\mathfrak{B}$ be the ordered basis $(\alpha_{i},\alpha_{j})$ for 
$V_{i,j}$, and for any linear mapping $T: V_{i,j} \rightarrow 
V_{i,j}$ let $[T]_{\mathfrak{B}}$ be the matrix for $T$ relative to 
$\mathfrak{B}$.  
Then 
\[ [S_{i}|_{V_{i,j}}]_{\mathfrak{B}} = 
\left(\begin{array}{cc} -1 & -a_{ij}\\ 0 & 1\end{array}\right), 
[S_{j}|_{V_{i,j}}]_{\mathfrak{B}} = 
\left(\begin{array}{cc} 1 & 0\\ -a_{ji} & -1\end{array}\right), 
[S_{i}S_{j}|_{V_{i,j}}]_{\mathfrak{B}} = 
\left(\begin{array}{cc} a_{ij}a_{ji}-1 & a_{ij}\\ -a_{ji} 
& -1\end{array}\right)\] 
Analysis of the eigenvalues for $X_{i,j} := 
[S_{i}S_{j}|_{V_{i,j}}]_{\mathfrak{B}}$ 
as in the proofs of Proposition 3.13 of \cite{Kac} and 
Proposition 1.3.21 of \cite{Kumar} shows that 
$X_{i,j}$ has infinite order when $a_{ij}a_{ji} 
\geq 4$, and hence $S_{i}S_{j}$ 
has infinite order as an element of $GL(V)$.  When $0 < a_{ij}a_{ji} < 4$, 
write $a_{ij}a_{ji} = 4\cos^{2}\theta$ for $\theta  := \pi/m_{ij}$.  
In this case check that $X_{i,j}$ has two distinct complex eigenvalues 
($e^{2i\theta}$ and $e^{-2i\theta}$).  It follows that $X_{i,j}$ has 
finite order $m_{ij}$.  
When $a_{ij}a_{ji} = 0$, then $X_{i,j}  = 
\left(\begin{array}{cc} -1 & 0\\ 0 & -1\end{array}\right)$, 
which clearly has order $m_{ij} = 2$. 
Now assume $0 \leq a_{ij}a_{ji} < 4$, and set $V_{i,j}' := 
\{v \in V\, |\, B(\alpha_{i},v) = 0 = 
B(\alpha_{j},v)\}$.  
One can easily check 
that $V_{i,j} \cap V_{i,j}' = \{0\}$.  The facts that $\dim V_{i,j} = 2$, 
$\dim V_{i,j}' \geq n-2$, and $V_{i,j} \cap V_{i,j}' = \{0\}$ together 
imply that $\dim V_{i,j}' = n-2$ and $V = V_{i,j} \oplus V_{i,j}'$. 
Since $S_{i}S_{j}$ acts as the identity on 
$V_{i,j}'$, it follows that 
$S_{i}S_{j}$ has order $m_{ij}$ as an element of $GL(V)$.  

Then  there is a unique homomorphism 
$\sigma_{A}: W \rightarrow GL(V)$ for which $\sigma_{A}(s_{i}) = 
S_{i}$. 
With the dependence on $A$ understood, we set $\sigma := \sigma_{A}$. 
We now have $W$ acting on $V$, and for all $w \in W$ and $v \in V$ we 
write $w.v$ for $\sigma(w)(v)$. 
We call $\sigma$ a {\em geometric representation} 
of $W$.  
If $A$ is symmetric such 
that $a_{kl}a_{lk} \geq 4 \Rightarrow a_{kl} = a_{lk} = -2$ for all $k 
\not= l$, then $\sigma$ is the {\em standard} geometric 
representation. 
The {\em root system} for 
$\sigma$ 
is $\Phi := \Phi_{A} := \{w.\alpha_{i}\}_{i 
\in I_{n}, w \in W}$.  
For each $w \in W$,  
$\sigma(w)$ permutes $\Phi$, so $\sigma$ induces an 
action of $W$ on $\Phi$. Evidently,  
$\Phi = -\Phi$.  Elements of $\Phi$ are 
{\em roots} and are necessarily nonzero.  
If $\alpha = \sum c_{i}\alpha_{i}$ is a root with all 
$c_{i}$ nonnegative (respectively nonpositive), then say $\alpha$ is a {\em 
positive} (resp.\ {\em negative}) root. 
Let $\Phi^{+}$ and $\Phi^{-}$ 
denote the collections of positive and negative roots respectively.  
Clearly $\Phi^{+} \cap \Phi^{-} = \emptyset$.  
The next statement is Proposition 4.2.5 of 
\cite{BB} and appears in a somewhat different form as Corollary 4.3 in 
\cite{ErikssonThesis}. 

\noindent 
{\bf \BBProposition}\ \ {\sl Let $w \in W$ and $i \in I_{n}$. If 
$\myl(ws_{i}) > 
\myl(w)$, then $w.\alpha_{i} \in \Phi^{+}$.  
If $\myl(ws_{i}) < 
\myl(w)$, then $w.\alpha_{i} \in \Phi^{-}$.} 

\noindent 
This result analogizes Theorem 5.4 of \cite{HumCoxeter}, which 
handles the standard case. 
As with Corollary 5.4 of \cite{HumCoxeter}, it is a consequence of  
\BBProposition\ 
that the representation $\sigma$ is faithful. (See \cite{BB} 
Theorem 4.2.7.)   
It also follows that $\Phi = \Phi^{+} \cup \Phi^{-}$.  
This is Equation 
4.24 of \cite{BB}, which actually could have been derived at the end 
of Section 4.2 of that text.    

Kac--Moody Weyl groups are subsumed into this paper as 
follows: 
Let $A$ be a generalized Cartan matrix.  
We identify our simple roots $\{\alpha_{1},\ldots,\alpha_{n}\}$ with 
the simple roots in $\mathfrak{h}_{\mathbb{R}}^{*}$ 
of \cite{Kac}, which is the dual of a 
real vector space $\mathfrak{h}_{\mathbb{R}}$ of dimension $n+l$, where $l = 
\mbox{nullity}(A)$.  The simple ``coroots'' of \cite{Kac} are a linearly 
independent set $\{\alpha_{1}^{\vee},\ldots,\alpha_{n}^{\vee}\} \subset 
\mathfrak{h}_{\mathbb{R}}$  
for which $\alpha_{j}(\alpha_{i}^{\vee}) = a_{ij}$.  Now for $1 \leq 
i \leq n$, a mapping $R_{i}: \mathfrak{h}_{\mathbb{R}}^{*} \rightarrow 
\mathfrak{h}_{\mathbb{R}}^{*}$ 
is defined in \cite{Kac} 
by $R_{i}(v) = v - v(\alpha_{i}^{\vee})\alpha_{i}$.  The associated 
Kac--Moody Weyl group is the subgroup of $GL(\mathfrak{h}_{\mathbb{R}}^{*})$ 
generated by $\{R_{i}\}_{i=1}^{n}$.  If we identify our $V$ with 
$\mbox{span}_{\mathbb{R}}\{\alpha_{1},\ldots,\alpha_{n}\} \subseteq 
\mathfrak{h}_{\mathbb{R}}^{*}$ and restrict each $R_{i}$ to $V$, then 
the homomorphism $W \rightarrow GL(V)$ determined by $s_{i} \mapsto 
R_{i}|_{V}$ is the representation $\sigma$.  The real roots of 
Kac--Moody theory are the roots $\Phi \subset 
V$ 
obtained here from this geometric representation of $W$. 

{\bf \S 3 Root system results.} 
Asymmetry of the bilinear form leads to 
crucial differences with the symmetric case.  
Most notably, 
$\sigma(W)$ preserves the form $B$ if and only if $A$ is 
symmetric. From this fact  
for symmetric $A$ it readily 
follows that if $K\alpha_{x} \in \Phi$ for some $x \in I_{n}$ 
and real number $K$, 
then $K = \pm{1}$. (See equation 4.27 of \cite{BB}.)   
However, when $A$ is asymmetric 
sometimes $K\alpha_{x}$ is 
a root for $K \not= \pm{1}$, as can be seen in 
Exercise 4.9 of 
\cite{BB} and \UnitalExample\ 
\renewcommand{\thefootnote}{3}
below.\footnote{In Proposition 6.9 of \cite{ErikssonThesis} and in 
\cite{ErikssonDiscrete} just prior to Proposition 4.4, 
it is asserted that $s_{x}(\Phi^{+} \setminus 
\{\alpha_{x}\}) = \Phi^{+} \setminus 
\{\alpha_{x}\}$ for all $x \in I_{n}$. However, this will not be the case if 
$K\alpha_{x}$ is a root for some $K \not= \pm{1}$.  Only 
Theorem 6.9 of \cite{ErikssonThesis} and Proposition 4.4 of 
\cite{ErikssonDiscrete} are affected by this misstatement.  
(See \PositiveToNegativeRootsLemma\ below.)}   
To understand how such a $W$-action can generate scalar multiples of roots 
in $\Phi$, we first analyze how 
$s_{i}$ and $s_{j}$ act in tandem on $V_{i,j}$. Our next result 
strengthens Lemma 4.2.4 of \cite{BB} and provides a different proof.  
It also answers Exercise 4.6 of \cite{BB}. 

\noindent 
{\bf \TwoGeneratorAnalysis}\ \ 
{\sl Fix $i \not= j$ in $I_{n}$, and let $k$ be a positive integer.   
If $m_{ij} = \infty$, then 
$(s_{i}s_{j})^{k}.\alpha_{i} = 
a\alpha_{i} + b\alpha_{j}$ and 
$s_{j}(s_{i}s_{j})^{k}.\alpha_{i} = c\alpha_{i} + d\alpha_{j}$   
for positive coefficients $a$, $b$, $c$, and $d$.  
Now suppose $m_{ij} < \infty$. 
If $2k<m_{ij}$, then $(s_{i}s_{j})^{k}.\alpha_{i} = a\alpha_{i} + 
b\alpha_{j}$ with $a \geq 0$ and $b > 0$.  In this case, $a = 0$ 
if and only if $m_{ij}$ is odd 
and $k = (m_{ij}-1)/2$, and consequently  
$(s_{i}s_{j})^{k}.\alpha_{i} = 
\frac{-a_{ji}}{2\cos(\pi/m_{ij})}\alpha_{j}$.  
Similarly, if $2k<m_{ij}-1$, then 
$s_{j}(s_{i}s_{j})^{k}.\alpha_{i} = c\alpha_{i} + d\alpha_{j}$ with $c > 
0$ and $d \geq 0$.  In this case, $d = 0$ if and only if $m_{ij}$ is even 
and $k = (m_{ij}-2)/2$, and consequently 
$s_{j}(s_{i}s_{j})^{k}.\alpha_{i} = \alpha_{i}$.}

{\em Proof.} 
Let $\mathfrak{B}$ and $X_{i,j}$ be as above, and set 
$X_{i} := [S_{i}|_{V_{i,j}}]_{\mathfrak{B}}$  
and $X_{j} := 
[S_{j}|_{V_{i,j}}]_{\mathfrak{B}}$.  
To understand $(s_{i}s_{j})^{k}.\alpha_{i}$ and 
$s_{j}(s_{i}s_{j})^{k}.\alpha_{i}$ we compute 
$X_{i,j}^{k}$ and $X_{j}X_{i,j}^{k}$. Set 
$p := -a_{ij}$ and $q := -a_{ji}$. 

For $m_{ij} = \infty$, first take $pq = 4$.  We can write $X_{i,j} = 
PYP^{-1}$ for nonsingular $P$ and upper triangular $Y$ as follows: 
\[X_{i,j} = \frac{1}{p}
\left(\begin{array}{cc}p & p\\ 
2 & 1\end{array}\right)
\left(\begin{array}{cc} 1 & 1\\ 0 & 
1\end{array}\right) 
\left(\begin{array}{cc} -1 & p\\ 2 
& -p\end{array}\right).\] 
Then for any positive integer $k$ we obtain 
$X_{i,j}^{k} = 
\left(\begin{array}{cc}2k+1 & -kp\\ 
kq & -2k+1\end{array}\right)$.  It follows that 
$(s_{i}s_{j})^{k}.\alpha_{i} = (2k+1)\alpha_{i} +  
kq\alpha_{j}$, with both coefficients of the linear combination 
positive. From the first column of the matrix $X_{j}X_{i,j}^{k}$ we 
see that  
$s_{j}(s_{i}s_{j})^{k}.\alpha_{i} = (2k+1)\alpha_{i} +  
(k+1)q\alpha_{j}$, with both coefficients of the linear combination 
positive.    
Next take $pq > 4$.  In this case we get distinct eigenvalues 
$\lambda = \frac{1}{2}(pq-2+\sqrt{pq(pq-4)}) > 1$ and  
$\mu = \frac{1}{2}(pq-2-\sqrt{pq(pq-4)}) < 1$ for $X_{i,j}$ (here we 
have $\lambda\mu = 1$).  Similar to the above, we may write
$X_{i,j} = PDP^{-1}$ for the diagonal 
matrix $D = \left(\begin{array}{cc}\lambda & 0\\ 
0 & \mu\end{array}\right)$ 
and a nonsingular matrix $P$, from which we obtain 
\[X_{i,j}^{k} = \frac{1}{p(\lambda-\mu)}
\left(\begin{array}{cc}p & p\\ 
\mu' & \lambda'\end{array}\right)
\left(\begin{array}{cc} \lambda^{k} & 0\\ 0 & 
\mu^{k}\end{array}\right) 
\left(\begin{array}{cc} \lambda' & -p\\ -\mu' 
& p\end{array}\right),\]
for any positive integer $k$, where $\lambda' := \lambda+1$  
and $\mu' := \mu+1$. This (eventually) simplifies to 
\[X_{i,j}^{k} =  
\frac{1}{\lambda-\mu}
\left(\begin{array}{cc}\lambda'\lambda^{k} - \mu'\mu^{k} & 
-p(\lambda^{k}-\mu^{k})\\ 
q(\lambda^{k}-\mu^{k}) & \lambda'\mu^{k} - 
\mu'\lambda^{k}\end{array}\right).\] From this we also get 
\[X_{j}X_{i,j}^{k} =  
\frac{1}{\lambda-\mu}
\left(\begin{array}{cc}\lambda'\lambda^{k} - \mu'\mu^{k} & 
-p(\lambda^{k}-\mu^{k})\\ 
q(\lambda^{k+1}-\mu^{k+1}) & \mu'\mu^{k} - 
\lambda'\lambda^{k}\end{array}\right).\] 
The factor $\frac{1}{\lambda-\mu}$ is positive, 
and for both matrices $X_{i,j}^{k}$ and 
$X_{j}X_{i,j}^{k}$, 
the first column entries are positive.  
So,  
$(s_{i}s_{j})^{k}.\alpha_{i} = a\alpha_{i} + b\alpha_{j}$ with both 
$a$ and $b$ positive, and 
$s_{j}(s_{i}s_{j})^{k}.\alpha_{i} = c\alpha_{i} + d\alpha_{j}$ with $c$ 
and $d$ both positive.

For the $m_{ij} < \infty$ case, set $\theta := 
\pi/m_{ij}$.  Note that the hypotheses of the lemma require that 
$m_{ij} > 2$, so in particular $p$ and $q$ are nonzero.   
Check that $X_{i,j}$ can be written as 
$X_{i,j} = PDP^{-1}$ for a 
nonsingular matrix $P$ and diagonal matrix $D$ in the following way: 
\[\frac{1}{q(e^{2i\theta} - e^{-2i\theta})}
\left(\begin{array}{cc}e^{2i\theta}+1 & e^{-2i\theta}+1\\ 
q & q\end{array}\right)
\left(\begin{array}{cc} e^{2i\theta} & 0\\ 0 & 
e^{-2i\theta}\end{array}\right) 
\left(\begin{array}{cc} q & -e^{-2i\theta}-1\\ -q 
& e^{2i\theta}+1\end{array}\right).\]
Then for any positive integer $k$ we have  
\[X_{i,j}^{k} = PD^{k}P^{-1} = \frac{1}{\sin(2\theta)}
\left(\begin{array}{cc}
\sin(2(k+1)\theta) + \sin(2k\theta) & 
-p\sin(2k\theta)\\
q\sin(2k\theta) & -\sin(2k\theta) - \sin(2(k-1)\theta)
\end{array}\right)\]
and 
\[X_{j}X_{i,j}^{k} = \frac{1}{\sin(2\theta)}
\left(\begin{array}{cc}
\sin(2(k+1)\theta) + \sin(2k\theta) & 
-p\sin(2k\theta)\\
q\sin(2(k+1)\theta) & (1-pq)\sin(2k\theta) + \sin(2(k-1)\theta)
\end{array}\right)
\] 
Use the first column of $X_{i,j}^{k}$ and 
$X_{j}X_{i,j}^{k}$ 
to see that 
$(s_{i}s_{j})^{k}.\alpha_{i} = \frac{1}{\sin(2\theta)}
[\sin(2(k+1)\theta) + 
\sin(2k\theta)]\alpha_{i} +  
\frac{q}{\sin(2\theta)}\sin(2k\theta)\alpha_{j}$ and that 
$s_{j}(s_{i}s_{j})^{k}.\alpha_{i} =  
\frac{1}{\sin(2\theta)}[\sin(2(k+1)\theta) + 
\sin(2k\theta)]\alpha_{i} +  
\frac{q}{\sin(2\theta)}\sin(2(k+1)\theta)\alpha_{j}$.  As long as 
$2(k+1) < m_{ij}$, then all the coefficients of these linear 
combinations will be positive.  So now suppose $2(k+1) \geq m_{ij}$.  
First we consider $(s_{i}s_{j})^{k} = a\alpha_{i} + b\alpha_{j}$ 
for some positive $k$ with $2k < 
m_{ij}$.  There are two possibilities now: $2(k+1) = m_{ij}$ or 
$2(k+1) = m_{ij}+1$.  In the former case 
both $a$ and $b$ are positive.  In the latter case we have $m_{ij}$ 
odd, $a = \frac{1}{\sin(2\theta)}
[\sin(2(k+1)\theta) + 
\sin(2k\theta)] = 0$, and $b = \frac{q\sin\theta}{\sin(2\theta)} = 
\frac{q}{2\cos\theta}$.  Second we consider 
$s_{j}(s_{i}s_{j})^{k} = c\alpha_{i} + d\alpha_{j}$ 
for some positive $k$ with $2k < m_{ij}-1$. Now the fact that $2(k+1) 
\geq m_{ij}$ implies we have $2(k+1) = m_{ij}$.  In particular, 
$m_{ij}$ is even.  With $k = (m_{ij}-2)/2$ now, one can check that 
$d = 0$ and $c = 1$.\hfill\QED

Distinct nodes $\gamma_{i}$ and $\gamma_{j}$ 
in $(\Gamma,A)$ are {\em odd-neighborly} if $m_{ij}$ is odd. 
If in addition we have 
$a_{ij} \not= a_{ji}$, then the adjacent nodes $\gamma_{i}$ 
and $\gamma_{j}$ form an {\em odd asymmetry}. 
For odd $m_{ij}$,   
let $v_{ji}$ be the element 
$(s_{i}s_{j})^{(m_{ij}-1)/2}$ of $W$, and 
set $K_{ji} := \frac{-a_{ji}}{2\cos(\pi/m_{ij})} = 
\sqrt{\frac{a_{ji}}{a_{ij}}}$. 
In view of \TwoGeneratorAnalysis,  
$v_{ji}.\alpha_{i} = K_{ji}\alpha_{j}$.  
Observe that $K_{ij}K_{ji} = 1$ and that 
$v_{ij} = v_{ji}^{-1}$.  
We have $\myl(v_{ji}) = m_{ij}-1$. 
Say a sequence $\mathcal{P} := 
[\gamma_{i_{0}},\gamma_{i_{1}},\ldots,\gamma_{i_{p}}]$ of nodes from 
$\Gamma$ is a {\em path of odd neighbors}, or {\em ON-path}, 
if consecutive nodes of $\mathcal{P}$ are odd neighbors.  
The ON-path $\mathcal{P}$ has length $p$, and we allow 
ON-paths to have length zero.  
We say $\gamma_{i_{0}}$ and 
$\gamma_{i_{p}}$ are the {\em start} and {\em end} nodes of the 
ON-path, 
respectively. 
Let $w_{_{\mathcal{P}}} \in W$ be the Coxeter group element 
$v_{i_{p}i_{p-1}}\cdots{v}_{i_{2}i_{1}}v_{i_{1}i_{0}}$, and 
let $\Pi_{_{\mathcal{P}}} := 
K_{i_{p}i_{p-1}}\cdots{K}_{i_{2}i_{1}}K_{i_{1}i_{0}}$, where 
$w_{_{\mathcal{P}}} = \varepsilon$ with $\Pi_{_{\mathcal{P}}} = 1$ 
when $\mathcal{P}$ has length zero. 
Then $w_{_{\mathcal{P}}}.\alpha_{i_{0}} = 
\Pi_{_{\mathcal{P}}}\alpha_{i_{p}}$.  If ON-path $\mathcal{Q} = 
[\gamma_{j_{0}},\gamma_{j_{1}},\ldots,\gamma_{j_{q}}]$ 
has the same start node as the end node of  
$\mathcal{P}$, then their {\em concatenation} 
$\mathcal{P}\natural\mathcal{Q}$ is the ON-path  
$[\gamma_{i_{0}},\gamma_{j_{1}},\ldots,
\gamma_{i_{p}}=\gamma_{j_{0}},\ldots,\gamma_{j_{q}}]$. 
Note that 
$w_{_{\mathcal{P}\natural\mathcal{Q}}} = 
w_{_{\mathcal{Q}}}w_{_{\mathcal{P}}}$. 

Distinct nodes $\gamma_{i}$ and $\gamma_{j}$ 
in $(\Gamma,A)$ are {\em 
even-related} if $m_{ij}$ is even. For even $m_{ij}$, let $v_{ji}$ be 
the element $s_{j}(s_{i}s_{j})^{(m_{ij}-2)/2}$ of $W$.  Then 
$v_{ji}.\alpha_{i} = \alpha_{i}$ (for $m_{ij} \geq 
4$ this is justified by \TwoGeneratorAnalysis), and $v_{ij} = v_{ji}^{-1}$. 
We have $\myl(v_{ji}) = m_{ij}-1$. 
Say a sequence $\mathcal{S} := [(\gamma_{i_{0}},\gamma_{i_{1}}), 
(\gamma_{i_{0}},\gamma_{i_{2}}), \ldots , (\gamma_{i_{0}},\gamma_{i_{p}})]$ 
is a {\em sequence of 
even-related nodes}, or {\em ER-sequence}, {\em rooted at} 
$\gamma_{i_{0}}$ if for each pair 
$(\gamma_{i_{0}},\gamma_{i_{k}})$ of the sequence $(1 \leq k \leq p)$ 
the nodes $\gamma_{i_{0}}$ and $\gamma_{i_{k}}$ are even-related.  
Say $\mathcal{S}$ has length $p$.  We allow $\mathcal{S}$ to be the 
empty sequence, in which case it has length zero.  
Let $w_{_{\mathcal{S}}} \in W$ be the Coxeter group element 
$v_{i_{p}i_{0}}\cdots{v}_{i_{2}i_{0}}v_{i_{1}i_{0}}$, with 
$w_{_{\mathcal{S}}} = \varepsilon$ 
when $\mathcal{S}$ has length zero. 
Then $w_{_{\mathcal{S}}}.\alpha_{i_{0}} = \alpha_{i_{0}}$. 
If ER-sequence $\mathcal{T} = 
[(\gamma_{i_{0}},\gamma_{j_{1}}), \ldots , 
(\gamma_{i_{0}},\gamma_{j_{q}})]$ is also rooted at $\gamma_{i_{0}}$, 
then the {\em concatenation} 
$\mathcal{S}\sharp\mathcal{T}$ is the ER-sequence  
$[(\gamma_{i_{0}},\gamma_{i_{1}}), \ldots , 
(\gamma_{i_{0}},\gamma_{i_{p}}), (\gamma_{i_{0}},\gamma_{j_{1}}), 
\ldots , (\gamma_{i_{0}},\gamma_{j_{q}})]$ rooted at $\gamma_{i_{0}}$. 
Note that 
$w_{_{\mathcal{S}\sharp\mathcal{T}}} = 
w_{_{\mathcal{T}}}w_{_{\mathcal{S}}}$. 

\noindent 
{\bf \HumphreysTheorem}\ \ 
{\sl Let $w \in W$ with $w \not= \varepsilon$, and let $i \in I_{n}$. 
(1) 
Then $w.\alpha_{i} = 
K\alpha_{x}$ for some $x \in I_{n}$ and some $K > 0$  
if and only if there is an ON-path $\mathcal{P} = 
[\gamma_{i_{0}=i},
\gamma_{i_{1}},\ldots,\gamma_{i_{p-1}},\gamma_{i_{p}=x}]$ and 
ER-sequences $\mathcal{S}_{k}$ rooted at $\gamma_{i_{k}}$  
$(0 \leq k \leq p)$ such that $w = 
w_{_{\mathcal{S}_{p}}} v_{i_{p},i_{p-1}} w_{_{\mathcal{S}_{p-1}}} 
v_{i_{p-1},i_{p-2}} \cdots v_{i_{2},i_{1}} w_{_{\mathcal{S}_{1}}} 
v_{i_{1},i_{0}} w_{_{\mathcal{S}_{0}}}$.  In this case, 
$w.\alpha_{i} = 
w_{_{\mathcal{P}}}.\alpha_{i} = \Pi_{_{\mathcal{P}}}\alpha_{x}$.  
(2) 
Similarly $w.\alpha_{i} = 
K\alpha_{x}$ for some $x \in I_{n}$ and some $K<0$  
if and only if 
there is an ON-path $\mathcal{P} = 
[\gamma_{i_{0}=i},
\gamma_{i_{1}},\ldots,\gamma_{i_{p-1}},\gamma_{i_{p}=x}]$ and 
ER-sequences $\mathcal{S}_{k}$ rooted at $\gamma_{i_{k}}$  
$(0 \leq k \leq p)$ such that $w = 
w_{_{\mathcal{S}_{p}}} v_{i_{p},i_{p-1}} w_{_{\mathcal{S}_{p-1}}} 
v_{i_{p-1},i_{p-2}} \cdots v_{i_{2},i_{1}} w_{_{\mathcal{S}_{1}}} 
v_{i_{1},i_{0}} w_{_{\mathcal{S}_{0}}} s_{i}$. In this case, 
$w.\alpha_{i} = (w_{_{\mathcal{P}}}s_{i}).\alpha_{i} = 
-\Pi_{_{\mathcal{P}}}\alpha_{x}$.} 

{\em Proof.} 
Note that (2) follows from (1).  For (1), the ``if'' 
direction is handled by the two definitions paragraphs preceding the 
theorem statement.  
For the ``only if'' direction, we induct on $\myl(w)$.  
If $\myl(w) = 1$, then it is clear that $w = s_{j}$ for some $j \not= i$ 
in $I_{n}$ and with $m_{ij} = 2$. So, $v_{ji} = s_{j}$. 
Taking $\mathcal{S}_{0} = 
[(\gamma_{i},\gamma_{j})]$, $\mathcal{P} = [\gamma_{i}]$, and 
$\mathcal{S}_{1}$ the empty sequence, then $w$ 
has the desired form. 
Now suppose $\myl(w) > 1$.  
Take any $j \in I_{n}$ 
for which $\myl(ws_{j}) = \myl(w) - 1$. 
Since $\myl(ws_{i}) > \myl(w)$, then $i 
\not= j$.   
Let $J := \{i,j\}$, and let $v^{J}$ be the 
unique element in $W^{J}$ and $v_{_{J}}$ the unique element in $W_{J}$ 
for which $w = v^{J}v_{_{J}}$.  
Then $\myl(w) = \myl(v^{J}) + \myl(v_{_{J}})$ by Proposition 2.4.4 
of \cite{BB}.  
Write $v_{_{J}}.\alpha_{i} = 
a\alpha_{i} + b\alpha_{j}$.  Since $\myl(ws_{i}) > \myl(w)$, then 
$\myl(v_{_{J}}s_{i}) > \myl(v_{_{J}})$, and hence 
$v_{_{J}}.\alpha_{i} \in \Phi^{+}$ (\BBProposition). 
So $a \geq 0$ and $b \geq 0$. 
Suppose $a > 0$ and $b > 0$.  Now $v^{J} \in W^{J}$ implies that 
$\myl(v^{J}s_{i}) > \myl(v^{J})$ and 
$\myl(v^{J}s_{j}) > \myl(v^{J})$, and hence 
$v^{J}.\alpha_{i} \in \Phi^{+}$ and $v^{J}.\alpha_{j} \in 
\Phi^{+}$ (\BBProposition). 
Write $v^{J}.\alpha_{i} = \sum_{y\in{I}_{n}}c_{y}\alpha_{y}$ ($c_{y} 
\geq 0$) and 
$v^{J}.\alpha_{j} = \sum_{y\in{I}_{n}}d_{y}\alpha_{y}$ ($d_{y} \geq 
0$).  Then $K\alpha_{x} = w.\alpha_{i} = 
v^{J}.(a\alpha_{i}+b\alpha_{j}) = 
\sum_{y\in{I}_{n}}(ac_{y}+bd_{y})\alpha_{y}$ implies that 
for all $y \not= x$,
$ac_{y}+bd_{y} = 0$ and hence $c_{y} = d_{y} = 0$.  
Then $v^{J}.\alpha_{i}$ and $v^{J}.\alpha_{j}$ are both  
multiples of $\alpha_{x}$.  But then $(v^{J})^{-1}.\alpha_{x}$ is a 
scalar multiple of $\alpha_{i}$ and of $\alpha_{j}$, which is absurd.  
So we must have $a=0$ or $b=0$. 
Then by \TwoGeneratorAnalysis, it follows that $m_{ij}$ is finite and 
$v_{_{J}} = v_{ji}$.  

If $m_{ij}$ is even, then $v_{ji}.\alpha_{i} = \alpha_{i}$.  So 
$v^{J}.\alpha_{i} = K\alpha_{x}$.  If $v^{J} = \varepsilon$, then 
take $\mathcal{S}_{0} = 
[(\gamma_{i},\gamma_{j})]$, $\mathcal{P} = [\gamma_{i}]$, and 
$\mathcal{S}_{1}$ the empty sequence to see that 
$w=v_{ji}$ has the desired form. 
Otherwise, since $\myl(v^{J}) < \myl(w)$ we may 
apply the induction hypothesis to $v^{J}$ to see that there 
is an ON-path $\mathcal{P} = 
[\gamma_{i_{0}=i},
\gamma_{i_{1}},\ldots,\gamma_{i_{p-1}},\gamma_{i_{p}=x}]$ and 
ER-sequences $\mathcal{S}_{k}$ rooted at $\gamma_{i_{k}}$  
$(0 \leq k \leq p)$ such that $v^{J} = 
w_{_{\mathcal{S}_{p}}} v_{i_{p},i_{p-1}} w_{_{\mathcal{S}_{p-1}}} 
v_{i_{p-1},i_{p-2}} \cdots v_{i_{2},i_{1}} w_{_{\mathcal{S}_{1}}} 
v_{i_{1},i_{0}} w_{_{\mathcal{S}_{0}}}$.  Let $\mathcal{S}'_{0} := 
[(\gamma_{i=i_{0}},\gamma_{j})]\sharp\mathcal{S}_{0}$. Then we get $w = 
w_{_{\mathcal{S}_{p}}} v_{i_{p},i_{p-1}} w_{_{\mathcal{S}_{p-1}}} 
v_{i_{p-1},i_{p-2}} \cdots v_{i_{2},i_{1}} w_{_{\mathcal{S}_{1}}} 
v_{i_{1},i_{0}} w_{_{\mathcal{S}'_{0}}}$, which has the desired form.  
On the other 
hand, if $m_{ij}$ is odd, then $v_{ji}.\alpha_{i} = 
K_{ji}\alpha_{j}$.  So $v^{J}.\alpha_{j} = 
\frac{K}{K_{ji}}\alpha_{x}$. If $v^{J} = \varepsilon$, then take 
$\mathcal{P} := [\gamma_{i},\gamma_{j}]$ with 
$\mathcal{S}_{0}$ and $\mathcal{S}_{1}$ empty ER-sequences to see 
that $w=v_{ji}$ has the desired form.  Otherwise apply the induction 
hypothesis to $v^{J}$ to see that there  
is an ON-path $\mathcal{P} = 
[\gamma_{i_{1}=j},
\gamma_{i_{2}},\ldots,\gamma_{i_{p-1}},\gamma_{i_{p}=x}]$ and 
ER-sequences $\mathcal{S}_{k}$ rooted at $\gamma_{i_{k}}$ 
$(1 \leq k \leq p)$ such that $v^{J} = 
w_{_{\mathcal{S}_{p}}} v_{i_{p},i_{p-1}} w_{_{\mathcal{S}_{p-1}}} 
v_{i_{p-1},i_{p-2}} \cdots v_{i_{2},i_{1}} w_{_{\mathcal{S}_{1}}}$. 
Take $i_{0}=i$, $\mathcal{P}' = 
[\gamma_{i},\gamma_{j}]\natural\mathcal{P}$ (an ON-path), and 
$\mathcal{S}_{0}$ an empty ER-sequence.  Then we get $w = 
w_{_{\mathcal{S}_{p}}} v_{i_{p},i_{p-1}} w_{_{\mathcal{S}_{p-1}}} 
v_{i_{p-1},i_{p-2}} \cdots v_{i_{2},i_{1}} w_{_{\mathcal{S}_{1}}} 
v_{i_{1},i_{0}} w_{_{\mathcal{S}_{0}}}$ as desired. Whether $m_{ij}$ 
is even or odd, we now see that $w.\alpha_{i} = 
w_{_{\mathcal{P}}}.\alpha_{i} = 
\Pi_{_{\mathcal{P}}}\alpha_{x}$.\hfill\QED 

The induction argument of the preceding proof can be viewed as a 
constructive method for obtaining the expression of the theorem 
statement for the Coxeter group element $w$. A further consequence of 
the proof is the following result about the length of $w$. It says, in 
effect, that if we write $w$ as a product of $v_{ji}$'s as prescribed 
in the theorem statement and then in 
turn write each such $v_{ji}$ as a shortest product of generators, 
the resulting expression for $w$ is reduced.   

\noindent 
{\bf \HumphreysCorollary} {\sl Suppose $w \in W$, $i,x \in I_{n}$, 
and $w.\alpha_{i} = K\alpha_{x}$ for some $K>0$.  Suppose  
$w = 
w_{_{\mathcal{S}_{p}}} v_{i_{p},i_{p-1}} w_{_{\mathcal{S}_{p-1}}} 
v_{i_{p-1},i_{p-2}} \cdots v_{i_{2},i_{1}} w_{_{\mathcal{S}_{1}}} 
v_{i_{1},i_{0}} w_{_{\mathcal{S}_{0}}}$ for an ON-path $\mathcal{P} = 
[\gamma_{i_{0}=i},
\gamma_{i_{1}},\ldots,\gamma_{i_{p-1}},\gamma_{i_{p}=x}]$ and 
ER-sequences $\mathcal{S}_{k}$ rooted at $\gamma_{i_{k}}$ 
$(0 \leq k \leq p)$  obtained by the method of the preceding 
proof.  For all $j,l \in I_{n}$, let $c(j,l)$ count the total number of 
occurrences of $(\gamma_{j},\gamma_{l})$ as consecutive nodes (in this 
order) in the 
ON-path $\mathcal{P}$ or as a pair in the ER-sequences $\mathcal{S}_{k}$ $(0 \leq k 
\leq p)$.  Then $\myl(w) = \sum_{j,l \in I_{n}}c(j,l)(m_{jl}-1)$.}\hfill\QED

\begin{figure}[t]
\begin{center}
\UnitalExampleFigure: A unital ON-cyclic E-GCM graph for 
\UnitalExample.\\
{\small The notation $\!\!\!$\CircleIntegerm$\!\!$ (resp.\ 
$\!\!\!\!$\CircleInfty$\!\!\!$) 
on an edge 
$\!\!\!\!\!\!\!\!\!$\TwoCitiesGraphWithoutLabels$\!\!\!\!\!\!\!\!\!$ 
indicates that $pq = 4\cos^{2}(\pi/m)$ (resp.\ $pq \geq 4$).} 

\ 

\vspace*{-0.1in}
\BowTie   
\end{center}

\vspace*{-0.35in}
\end{figure}

For any $\alpha \in \Phi$, 
set $\mathfrak{S}(\alpha) := \mathfrak{S}_{A}(\alpha) := 
\{K\alpha\}_{K\in\mathbb{R}}\cap\Phi^{+}$. 
Our analysis of the sets 
$\mathfrak{S}(\alpha)$ requires some 
additional notation.  
For ON-paths $\mathcal{P}$ and $\mathcal{Q}$, 
write $\mathcal{P} \sim 
\mathcal{Q}$ and say $\mathcal{P}$ and $\mathcal{Q}$ are $\Pi$-{\em 
equivalent}  
if these ON-paths have 
the same start and end nodes and 
$\Pi_{_{\mathcal{P}}} = \Pi_{_{\mathcal{Q}}}$.  This is an 
equivalence relation on the set of all ON-paths. 
An ON-path $\mathcal{P}$ 
is {\em simple} if it has no repeated nodes with the possible 
exception that the start and end nodes may coincide.  
Two ON-paths $\mathcal{P}$ and $\mathcal{Q}$ are 
{\em scalar-distinct} if $\Pi_{_{\mathcal{P}}} 
\not= \Pi_{_{\mathcal{Q}}}$. 
An ON-path 
$\mathcal{P} = [\gamma_{i_{0}},\ldots,\gamma_{i_{p}}]$ 
is an {\em ON-cycle} if $\gamma_{i_{p}} = \gamma_{i_{0}}$.  
It is {\em unital} if 
$\Pi_{\mathcal{P}} = 1$, i.e.\ 
$a_{i_{0},i_{1}}a_{i_{1},i_{2}}{\cdots}a_{i_{p-1},i_{0}} = 
a_{i_{0},i_{p-1}}{\cdots}a_{i_{2},i_{1}}{\cdots}a_{i_{1},i_{0}}$.  
We say 
$(\Gamma,A)$ is {\em unital ON-cyclic} if and only if 
$\Pi_{_{\mathcal{C}}} = 1$ for all ON-cycles $\mathcal{C}$.  
See \UnitalExampleFigure. 
From the definitions it follows that  
$(\Gamma,A)$ is unital ON-cyclic if it has no odd asymmetries.  So if 
$A$ is a GCM, then $(\Gamma,A)$ is unital ON-cyclic. 
If $A$ is a
symmetrizable E-GCM, then by applying Exercise 2.1 of \cite{Kac} or Exercise 1.5.E.1 
of \cite{Kumar} to the environment of E-GCM's, one sees that 
$(\Gamma,A)$ is 
unital ON-cyclic. 
However, a unital ON-cyclic E-GCM graph need not have a symmetrizable 
matrix $A$, as \UnitalExample\ shows.  
To check if an E-GCM graph is unital ON-cyclic, it is enough to check 
that each simple ON-cycle is unital. 
An E-GCM graph is {\em ON-connected} if 
any two nodes can be joined by an ON-path.  An {\em ON-connected 
component} of 
$(\Gamma,A)$ is an E-GCM subgraph 
$(\Gamma',A')$ whose nodes 
form a maximal collection of nodes in $(\Gamma,A)$ which 
can be pairwise joined by ON-paths.  

\noindent 
{\bf \OneToOneLemma}\ \ {\sl Let 
$\alpha$ and $\beta$ be roots in $\Phi$.  Suppose an element of 
$\mathfrak{S}(\alpha)$ is in the same orbit as an element of 
$\mathfrak{S}(\beta)$ under 
the action of $W$ on $\Phi$.  
Then there is a one-to-one correspondence between the sets 
$\mathfrak{S}(\alpha)$ and $\mathfrak{S}(\beta)$.  
If $\gamma_{i}$ and $\gamma_{j}$ are nodes in the 
same ON-connected component of $(\Gamma,A)$,   
then there is a one-to-one 
correspondence between the sets 
$\mathfrak{S}(\alpha_{i})$ and 
$\mathfrak{S}(\alpha_{j})$.} 

{\em Proof.} Since $\mathfrak{S}(\alpha) = 
\mathfrak{S}(K\alpha)$ for all $K\alpha \in 
\mathfrak{S}(\alpha)$, it suffices to assume that $\alpha$ 
and $\beta$ are in the same $W$-orbit, i.e.\ $\beta = w.\alpha$ for 
some $w \in W$.  It is easy to see that the 
mapping $\mathfrak{S}(\alpha) \rightarrow 
\mathfrak{S}(\beta)$ given by 
$\sigma(w)|_{\mathfrak{S}(\alpha)}$ gives the desired 
one-to-one correspondence.  If $\gamma_{i}$ and $\gamma_{j}$ are in 
the same ON-connected component, then by \HumphreysTheorem, 
some positive scalar multiple 
of $\alpha_{j}$ is in the $W$-orbit of $\alpha_{i}$.  
Thus there is a a one-to-one 
correspondence between the sets 
$\mathfrak{S}(\alpha_{i})$ and 
$\mathfrak{S}(\alpha_{j})$.\hfill\QED

The proof of the following lemma is a routine verification, so it is omitted. 

\noindent 
{\bf \ReducedLemma}\ \ 
{\sl Suppose $(\Gamma,A)$ is unital ON-cyclic.  Then 
for any ON-path $\mathcal{P}$ there is a simple ON-path  
which is $\Pi$-equivalent to $\mathcal{P}$.}


Although \TFAE\ and \TFAECorollary\ ask readers to look at a 
subgraph $(\Gamma',A')$ of $(\Gamma,A)$, the conclusions pertain 
to the action of $W = W(\Gamma,A)$ on $\Phi$. 

\noindent 
{\bf \TFAE}\ \ 
{\sl Choose any ON-connected component $(\Gamma',A')$ of 
$(\Gamma,A)$, and let $J := \{x \in I_{n}\}_{\gamma_{x} \in \Gamma'}$. 
Then the following are equivalent:} 

\hspace{0.35in} {\sl (1) $(\Gamma',A')$ is unital ON-cyclic.} 

\hspace{0.35in} {\sl (2) $|\mathfrak{S}(w.\alpha_{x})| < 
\infty$ for some $x \in J$ and $w \in W$.} 

\hspace{0.35in} {\sl (3) $|\mathfrak{S}(w.\alpha_{x})| < 
\infty$ for all $x \in J$ and $w \in W$.}   

\noindent 
{\sl In these cases for all $x, y \in 
J$ and $w \in W$,  
$|\mathfrak{S}(\alpha_{x})| = 
|\mathfrak{S}(w.\alpha_{y})|$. This common quantity 
is equal to the largest 
number of pairwise scalar-distinct 
simple ON-paths in $(\Gamma,A)$ with end node $\gamma_{x}$.}

{\em Proof.}  We show (2) $\Rightarrow$ (1) $\Rightarrow$ (3), the 
implication (3) $\Rightarrow$ (2) being obvious.  
For (1) 
$\Rightarrow$ (3), let $x \in J$.  
Observe that if $K\alpha_{x} \in \Phi^{+}$, then 
by \HumphreysTheorem\ we must have $K = \Pi_{_{\mathcal{P}}}$ for 
some ON-path $\mathcal{P}$ with end node $\gamma_{x}$.  Therefore 
$\mathcal{P}$ is in $(\Gamma',A')$.  
By \ReducedLemma, we may take a simple  
ON-path $\mathcal{Q}$ $\Pi$-equivalent to $\mathcal{P}$ 
(all ON-paths $\Pi$-equivalent to 
$\mathcal{P}$ must be in $(\Gamma',A')$), so that $K = 
\Pi_{_{\mathcal{Q}}}$.  Since there can be at most a finite number of 
simple ON-paths, then there can be at most finitely many positive roots that 
are scalar multiples of a given $\alpha_{x}$.  
That $|\mathfrak{S}(w.\alpha_{x})| = |\mathfrak{S}(\alpha_{x})|$ for 
all $w \in W$ follows from \OneToOneLemma. 
For (2) $\Rightarrow$ 
(1), we show the contrapositive.  Let $\mathcal{C} = 
[\gamma_{x},\ldots,\gamma_{x}]$ be a non-unital ON-cycle with 
start/end node $\gamma_{x}$ for an $x \in J$. So 
necessarily $\mathcal{C}$ has nonzero length.  Note that 
$w_{_{\mathcal{C}}}.\alpha_{x} = \Pi_{_{\mathcal{C}}}\alpha_{x}$. 
Next, for $y \in J$ (and possibly $y = 
x$) take any ON-path $\mathcal{P}$ 
with start node $\gamma_{x}$ 
and end node $\gamma_{y}$.  Since 
$w_{_{\mathcal{P}}}.\alpha_{x} = \Pi_{_{\mathcal{P}}}\alpha_{y}$, it 
follows that $w_{_{\mathcal{P}}}w_{_{\mathcal{C}}}^{k}.\alpha_{x} = 
\Pi_{_{\mathcal{P}}}\Pi_{_{\mathcal{C}}}^{k}\alpha_{y}$ for any integer $k$.  
In particular, for all $y \in J$, we have 
$|\mathfrak{S}(\alpha_{y})| = \infty$. 
So by \OneToOneLemma\ $|\mathfrak{S}(w.\alpha_{y})| = \infty$ for 
all $y \in J$, $w \in W$. 
The next-to-last 
claim of the theorem statement follows from \OneToOneLemma. The final 
claim follows from our proof above of the (1) $\Rightarrow$ (3) part 
of the theorem statement.\hfill\QED 

From \HumphreysTheorem\ it follows 
that {\sl if 
$(\Gamma,A)$ has 
an odd asymmetry, then there exists 
a root which is a 
non-trivial multiple of a simple root}.  The following corollary of 
\TFAE\ contains a more general statement 
that includes the converse. 
When $A$ is an integer matrix, odd neighbors $\gamma_{i}$ and 
$\gamma_{j}$ must have $\{-a_{ij},-a_{ji}\} = \{1,1\}$.  These are not 
asymmetric.  Therefore the matrices $A$ defining Weyl groups have no 
odd asymmetries. 
In this integer matrix setting, Kac (\cite{Kac} Proposition 5.1.b) and 
Kumar (\cite{Kumar} Corollary 1.3.6.a) show that for a ``real'' root 
$\alpha$ and real number $K$, $K\alpha$ is also a root if and only if 
$K = \pm{1}$.  Their proofs use Lie brackets and root space 
reasoning.  But alternatively, 
this result is also a very special case of the following: 

\noindent 
{\bf \TFAECorollary}\ \ {\sl We have 
$|\mathfrak{S}(\alpha)| = 1$ for all $\alpha \in \Phi$  
if and only if $(\Gamma,A)$ has no odd asymmetries. More generally, 
choose any ON-connected component $(\Gamma',A')$ of 
$(\Gamma,A)$, and let $J := \{x \in I_{n}\}_{\gamma_{x} \in \Gamma'}$. Then 
$|\mathfrak{S}(w.\alpha_{x})| = 1$ for some $x \in J$ and $w \in W$ 
if and only if 
$|\mathfrak{S}(w.\alpha_{x})| = 1$ for all $x \in J$ and $w \in W$ 
if and only if $(\Gamma',A')$ has no odd asymmetries.} 

{\em Proof.} Follows from \TFAECorollaryTheorems.\hfill\QED

Analogizing \cite{BB} and \cite{HRT}, for any $w \in W$ set $N(w) := 
N_{A}(w) := 
\{\alpha \in \Phi^{+}\, |\, w.\alpha \in \Phi^{-}\}$.  (For 
the matrices $A$ considered in \cite{Kumar} Ch.\ 1, this set   
is notated $\Phi_{w^{-1}}$.) 

\noindent
{\bf \PositiveToNegativeRootsLemma}\ \ {\sl For any $i \in I_{n}$, 
$s_{i}(\Phi^{+}\setminus\mathfrak{S}(\alpha_{i})) = 
\Phi^{+}\setminus\mathfrak{S}(\alpha_{i})$.  Now let $w \in 
W$.  
If $w.\alpha_{i} \in \Phi^{+}$, then $N(ws_{i}) = s_{i}(N(w)) 
\disjointunion 
\mathfrak{S}(\alpha_{i})$, a disjoint union. 
If $w.\alpha_{i} \in \Phi^{-}$, then $N(ws_{i}) = 
s_{i}(N(w)\setminus\mathfrak{S}(\alpha_{i}))$.} 

{\em Proof.}   Using \BBProposition, 
the proof of Proposition 5.6.(a) from \cite{HumCoxeter} 
is easily adjusted to prove the first claim.  Proofs for 
the remaining claims involve routine set inclusion 
arguments.\hfill\QED 

When $(\Gamma,A)$ is ON-connected and unital ON-cyclic, 
set 
$f_{\Gamma,A} := |\mathfrak{S}(\alpha_{x})|$ for 
any given $x \in I_{n}$. At this point, \BBProposition,   
\HumphreysTheorem, \PositiveToNegativeRootsLemma, and \TFAE\ 
allow us to modify the proof of Proposition 5.6 of \cite{HumCoxeter} 
to obtain the result that for all $w \in W$, 
$|N(w)| = f_{\Gamma,A}\, \myl(w)$. 
\LengthProposition\ below generalizes this statement for arbitrary 
E-GCM graphs.    
For $J \subseteq I_{n}$, 
let $\mathfrak{C}(J)$ denote the set of all 
ON-connected components of $(\Gamma,A)$ containing  
some node from the set $\{\gamma_{x}\}_{x \in J}$. 

\noindent 
{\bf \LengthProposition}\ \ 
{\sl Let $w \in W$ with $p = \myl(w) > 0$.  (1) Then $N(w)$ is 
finite if and only if $w$ has a reduced expression 
$s_{i_{1}}{\cdots}s_{i_{p}}$ for which 
$\mathfrak{S}(\alpha_{i_{q}})$ is finite for all $1 \leq q \leq p$ 
if and only if every 
reduced expression $s_{i_{1}}{\cdots}s_{i_{p}}$ for $w$ has 
$\mathfrak{S}(\alpha_{i_{q}})$ finite for all $1 \leq q \leq p$. 
(2) Now suppose $w = s_{i_{1}}{\cdots}s_{i_{p}}$ and 
$N(w)$ is finite.  Let $J := \{i_{1},\ldots,i_{p}\}$. 
In view of (1), let 
$f_{1}$ be the min and $f_{2}$ the max of all integers in the set 
$\{f_{\Gamma',A'}\, |\, (\Gamma',A') \in \mathfrak{C}(J)\}$.  
Then $f_{1}\, \myl(w) \leq |N(w)| \leq f_{2}\, \myl(w)$.} 
 
{\em Proof.} (1) follows from  
\PositiveToNegativeRootsLemma. For (2), 
induct on $\myl(w)$.  
Take $w' := s_{i_{1}}\cdots{s}_{i_{p-1}}$ 
with $w = w's_{i_{p}}$. 
Now $\gamma_{i_{p}}$ is 
in an ON-connected component $(\Gamma',A')$ 
of $(\Gamma,A)$.  Then by 
\PositiveToNegativeRootsLemma, $|N(w)| = |N(w')| + 
f_{\Gamma',A'}$.  Since $f_{1}\, \myl(w') \leq |N(w')| \leq 
f_{2}\, \myl(w')$, the result follows.\hfill\QED 

Apply \ListForCorollary\ to get: 

\noindent 
{\bf \UnitalCorollary}\ \ {\sl We have $N(w)$ finite for all $w \in 
W$ if and only if $(\Gamma,A)$ is unital 
ON-cyclic.  Moreover $|N(w)| = \myl(w)$ for all $w \in W$ if and only 
if $(\Gamma,A)$ has no odd asymmetries.}\hfill\QED

When $W$ is infinite, the length function must take arbitrarily large 
values. Then by \LengthProposition, 
$\Phi$ is infinite as well.  
If $W$ is finite, then $\Phi$ is finite as well, so 
$|\mathfrak{S}(\alpha_{x})| < \infty$ for all $x \in I_{n}$. 
In this case let $w_{0}$ 
be the longest element of $W$ (cf.\ Exercise 5.6.2 of \cite{HumCoxeter}). 
It is easily seen that if $w_{0} = s_{i_{1}}\cdots{s}_{i_{l}}$ is 
reduced, then $\{i_{1},\ldots,i_{l}\} = I_{n}$.  

\noindent 
{\bf \LengthCorollary}\ \ 
{\sl Suppose $W$ is finite. Let} $\Phi_{\mbox{\scriptsize std}}$ 
{\sl denote the root system for the 
standard geometric representation. 
Then} $f_{1}|\Phi_{\mbox{\scriptsize std}}^{+}| \leq 
|\Phi^{+}| \leq f_{2}|\Phi_{\mbox{\scriptsize std}}^{+}|$, 
{\sl where $f_{1}$ is the min and $f_{2}$ is the max of all integers in the set 
$\{f_{\Gamma',A'}\, |\, (\Gamma',A') \in \mathfrak{C}(I_{n})\}$.}

{\em Proof.} 
Apply \BBProposition\ to see that $N(w_{0}) = \Phi^{+}$.   By 
\LengthProposition, $f_{1}\myl(w_{0}) \leq 
|\Phi^{+}| \leq f_{2}\myl(w_{0})$.  To see that $\myl(w_{0}) = 
|\Phi_{\mbox{\scriptsize std}}^{+}|$, apply the previous reasoning in the standard 
case.\hfill\QED

\noindent 
{\bf \UnitalExample}\ \ In \UnitalExampleFigure\ is depicted a connected, 
unital ON-cyclic E-GCM graph $(\Gamma,A)$ 
with three ON-connected components: $(\Gamma_{1},A_{1})$ is the E-GCM 
subgraph with nodes $\gamma_{1}$ and $\gamma_{2}$;  
$(\Gamma_{2},A_{2})$ has nodes $\gamma_{4}$, $\gamma_{5}$, and 
$\gamma_{6}$; and $(\Gamma_{3},A_{3})$ has only the node $\gamma_{3}$.  
The matrix $A$ is not symmetrizable by Exercise 2.1 of 
\cite{Kac} or Exercise 1.5.E.1 
of \cite{Kumar}. 
Pertaining to the pair $(\gamma_{4},\gamma_{6})$, we have 
$4\cos^{2}(\pi/5) = \frac{3+\sqrt{5}}{2}$ and $2\cos(\pi/5) = 
\frac{1+\sqrt{5}}{2}$.  Since $a_{46} = -\frac{1+\sqrt{5}}{4}$ and 
$a_{64} = -(1+\sqrt{5})$, 
then $K_{46} = \frac{-a_{46}}{2\cos(\pi/5)} = 
\frac{1}{2}$ and $K_{64} = \frac{-a_{64}}{2\cos(\pi/5)} = 2$.  
For all other odd neighbors  
$(\gamma_{i},\gamma_{j})$ in this graph, $m_{ij} = 3$, so $K_{ij} = 
-a_{ij}$ and $K_{ji} = -a_{ji}$.  
By the last statement of  
\TFAE, 
$f_{\Gamma_{1},A_{1}} = 2$ and $f_{\Gamma_{2},A_{2}} = 3$.  For 
example, 
$\mathfrak{S}(\alpha_{2}) = 
\{\alpha_{2}, \frac{1}{5}\alpha_{2}\} = N(s_{2})$ and 
$\mathfrak{S}(\alpha_{5}) = 
\{\alpha_{5}, \frac{1}{7}\alpha_{5}, \frac{2}{7}\alpha_{5}\} = 
N(s_{5})$.  
By \LengthProposition, we can see that 
$f_{\Gamma_{1},A_{1}}\myl(s_{5}s_{2}) = 4 \leq |N(s_{5}s_{2})| 
\leq 6 = f_{\Gamma_{2},A_{2}}\myl(s_{5}s_{2})$. 
More precisely, by \PositiveToNegativeRootsLemma\ we get 
$N(s_{5}s_{2}) = s_{2}(N(s_{5}))
\disjointunion 
\mathfrak{S}(\alpha_{2})$,
whence $|N(s_{5}s_{2})| = 5$.\hfill\QED 

We now apply \HumphreysTheorem\ to extend a finiteness result of Brink and 
Howlett concerning a natural partial order on positive roots, cf.\ 
Theorem 2.8 of \cite{BH}.  
From here on, $\Phi_{\mbox{\scriptsize std}}$ denotes the root 
system for the standard geometric representation, and 
$\{\alpha_{i}^{\mbox{\scriptsize 
std}}\}_{i \in I_{n}}$ are its simple roots.   
Following \cite{BH} and \S 4.7 of 
\cite{BB}, for roots $\alpha, \beta \in \Phi_{\mbox{\scriptsize std}}^{+}$,  we say 
$\alpha$ {\em dominates} $\beta$, and write $\alpha\,\, 
\mathrm{dom}\,\,  
\beta$ if for all $w 
\in W$ we have $w.\beta \in \Phi_{\mbox{\scriptsize std}}^{-}$ whenever $w.\alpha \in 
\Phi_{\mbox{\scriptsize std}}^{-}$. It is known that 
the relation ``$\mathrm{dom}$'' on $\Phi_{\mbox{\scriptsize 
std}}^{+}$ 
is a partial order on $\Phi_{\mbox{\scriptsize std}}^{+}$ (\cite{BH}, 
\S 4.7 of \cite{BB}).  Roots in 
$\Phi_{\mbox{\scriptsize std}}^{+}$ that are minimal with respect to 
this partial order are {\em dominance-minimal}.  Observe that 
simple roots in $\Phi_{\mbox{\scriptsize std}}^{+}$ 
are dominance-minimal. The fact, due to Brink and Howlett 
in \cite{BH}, that the set 
of dominance-minimal elements is finite is viewed by some 
to be a fundamental result (see \cite{Casselman}, \S 4.7 of \cite{BB}).  Notable 
consequences of this finiteness result are the so-called Parallel Wall Theorem 
(discussed in \cite{BH}, see also \cite{Caprace}) for the associated 
Davis complex, as well as the fact that Coxeter groups are automatic 
\cite{BH}.  
However, the above 
definition of dominance does not extend nicely in the obvious way to the asymmetric 
setting: If adjacent nodes $\gamma_{i}$ and $\gamma_{j}$ in 
$(\Gamma,A)$ form an odd 
asymmetry, then by \HumphreysTheorem,  
$\alpha_{i}$ and $K\alpha_{i}$ are both positive 
roots for some positive $K\not=1$.  Then, each root would dominate the other, so  
dominance would not be an anti-symmetric relation on $\Phi^{+}$. 

In what follows, we address this issue.  
In the general setting,  
let $\Psi^{+} := \{\mathfrak{S}(\alpha)\}_{\alpha \in \Phi^{+}}$. 
For $\alpha, \beta \in \Phi^{+}$, say $\mathfrak{S}(\alpha)$ {\em dominates} 
$\mathfrak{S}(\beta)$ and write $\mathfrak{S}(\alpha)\,\,  
\mathrm{dom}\,\, \mathfrak{S}(\beta)$ if for all $w 
\in W$ we have $w.\beta \in \Phi^{-}$ whenever $w.\alpha \in 
\Phi^{-}$. It is easy to see that this definition is independent of 
the choice of representatives from each of $\mathfrak{S}(\alpha)$ and 
$\mathfrak{S}(\beta)$, so 
dominance is a well-defined 
relation on $\Psi^{+}$.    

\noindent 
{\bf \BrinkHowlettTheorem}\ \ {\sl Define a function} $\rho:\Psi^{+} 
\longrightarrow \Phi_{\mbox{\scriptsize 
std}}^{+}$ {\sl by the rule:} $\rho(\mathfrak{S}(\alpha)) := 
w.\alpha_{i}^{\mbox{\scriptsize std}}$ {\sl if $\alpha = w.\alpha_{i} 
\in \Phi^{+}$ for some $w \in W$ and $i \in I_{n}$.  Then $\rho$ 
is a well-defined 
bijection, and moreover  
$\rho$ and $\rho^{-1}$ both preserve dominance. In particular, 
``$\mathrm{dom}$'' is a partial order on $\Psi^{+}$ and the set 
of dominance-minimal elements of $\Psi^{+}$ 
is finite.}   

{\em Proof.} Obviously $\rho$ is 
surjective if it is well-defined.  
To see that $\rho$ is well-defined and injective, we show that for 
any $w_{1}, w_{2} \in W$ and  $i, x \in I_{n}$, we have 
$w_{1}.\alpha_{i} = w_{2}.\alpha_{x} \in \Phi^{+}$ if and only if 
$w_{1}.\alpha_{i}^{\mbox{\scriptsize std}} = w_{2}.\alpha_{x}^{\mbox{\scriptsize 
std}} \in \Phi_{\mbox{\scriptsize std}}^{+}$. Now if 
$w_{1}.\alpha_{i} = w_{2}.\alpha_{x} \in \Phi^{+}$, then 
$w.\alpha_{i} = \alpha_{x}$ for $w := (w_{2})^{-1}w_{1}$.  By 
\HumphreysTheorem, there is an ON-path $\mathcal{P} = 
[\gamma_{i_{0}=i},
\gamma_{i_{1}},\ldots,\gamma_{i_{p-1}},\gamma_{i_{p}=x}]$ and 
ER-sequences $\mathcal{S}_{k}$ rooted at $\gamma_{i_{k}}$ 
($0 \leq k \leq p$) such that $w = 
w_{_{\mathcal{S}_{p}}} v_{i_{p},i_{p-1}} w_{_{\mathcal{S}_{p-1}}} 
v_{i_{p-1},i_{p-2}} \cdots v_{i_{2},i_{1}} w_{_{\mathcal{S}_{1}}} 
v_{i_{1},i_{0}} w_{_{\mathcal{S}_{0}}}$. Using this expression for 
$w$, we can calculate that $w.\alpha_{i}^{\mbox{\scriptsize std}} = 
\alpha_{x}^{\mbox{\scriptsize std}}$. It follows that 
$w_{1}.\alpha_{i}^{\mbox{\scriptsize std}} = 
w_{2}.\alpha_{x}^{\mbox{\scriptsize 
std}}$, which is therefore positive in $\Phi_{\mbox{\scriptsize 
std}}$ by \BBProposition.  
The converse is entirely similar. Using \BBProposition, we have that 
$w.\alpha_{i} \in \Phi^{+}$ (respectively, $\Phi^{-}$) if and only if 
$w.\alpha_{i}^{\mbox{\scriptsize 
std}} \in \Phi_{\mbox{\scriptsize std}}^{+}$ (resp.\ 
$\Phi_{\mbox{\scriptsize std}}^{-}$). 
It now follows from the definitions 
that $\rho$ and $\rho^{-1}$ preserve dominance.  So, 
``$\mathrm{dom}$'' is a partial order on $\Psi^{+}$. 
That the set of 
dominance-minimal elements of $\Psi^{+}$ is finite now follows from 
Theorem 2.8 of \cite{BH}.\hfill\QED

{\bf \S 4 An application concerning the Tits cone.} 
We close with results which relate the size of a Coxeter group 
$W$ and 
the behavior of a ``fundamental domain'' for the ``contragredient''  
$W$-action.  The main result of this section (\TitsConeFiniteResult) 
is derived in two ways as an application of 
\UnitalCorollary/\HRTResult: first using the perspective of the 
numbers game, and second borrowing some results from \cite{Vinberg}. 
We continue to consider $\sigma: W \rightarrow GL(V)$. 
We have the natural pairing $\langle \lambda, v \rangle := 
\lambda(v)$ for elements $\lambda$ in the dual space $V^{*}$ and 
vectors $v$ in $V$. 
The contragredient representation 
$\sigma^{*} := \sigma_{A}^{*}: W \rightarrow GL(V^{*})$ 
is determined by $\langle \sigma^{*}(w)(\lambda), v \rangle = 
\langle \lambda, \sigma(w^{-1})(v) \rangle$.  When $w 
\in W$ and $\lambda \in V^{*}$, we 
write $w.\lambda$ for $\sigma^{*}(w)(\lambda)$.  
Let $D := \{\lambda \in V^{*}\, |\, \langle 
\lambda, \alpha_{i} \rangle \geq 0 \mbox{ for all } i \in I_{n}\}$.
Following \cite{Vinberg}, \cite{ErikssonThesis}, \cite{ErikssonDiscrete}, 
the {\em Tits 
cone} is $U := U_{A} := \cup_{w \in W}wD$. This generalizes the 
standard case of  
\cite{HumCoxeter}.  
In view of 
\BBProposition, the results of \cite{HumCoxeter} \S 5.13 hold 
here.  So $D$ is the aforementioned fundamental domain, and $U$ is a 
convex cone.  
Let $\overline{U}$ denote the closure of $U$. 
See the lecture notes of Howlett \cite{Howlett} for 
further discussion of properties of the Tits cone for the standard 
geometric representation $\sigma$, 
and in particular an investigation of phenomena in $\overline{U} 
\setminus U$.  
If $T$ is any convex cone, let $T_{0}$ denote the maximal subspace 
contained in $T$.  It is not hard to see that $T_{0} = T \cap (-T)$.  

Our Tits cone results below 
concern $U_{0}$.  These results both use/produce 
consequences from/for the numbers game. 
Elements of $V^{*}$ will be referred to as 
{\em positions} for 
$(\Gamma,A)$. 
We define a process of acting on positions in $V^{*}$ with certain 
sequences of Coxeter group generators that is equivalent to 
Eriksson's numbers game as presented in \S 4.3 of \cite{BB}. 
For a positive integer $p$ 
we say a sequence 
$(\gamma_{i_{1}},\ldots,\gamma_{i_{p}})$ from $(\Gamma,A)$ is {\em 
legal} from a given position $\lambda$ if 
$\langle s_{i_{q-1}}{\cdots}s_{i_{1}}.\lambda , \alpha_{i_{q}} 
\rangle > 0$ for all $1 \leq q \leq p$.  
Repeated application of \BBProposition\ 
implies that in this case, 
$s_{i_{p}}{\cdots}s_{i_{1}}$ is reduced. 
Call this 
the Reduced 
Word Result.  Next, say a position 
$\lambda$ is {\em good} if $\lambda \in -D$ or 
there exists a legal sequence 
$(\gamma_{i_{1}},\ldots,\gamma_{i_{p}})$ from $\lambda$ such that 
$s_{i_{p}}{\cdots}s_{i_{1}}.\lambda \in -D$.  In the latter 
case say $(\gamma_{i_{1}},\ldots,\gamma_{i_{p}})$ is a {\em 
terminated} legal sequence. 
Think of a good position as a position from which there is a 
(possibly empty) 
terminated legal sequence.  Eriksson's Strong Convergence Theorem 
(see Theorem 2.2 of \cite{ErikssonDiscrete}) 
shows that all legal sequences of maximal length 
from a good $\lambda$ terminate at the same ``terminal position'' in the same 
finite number of steps. 
Lemma 5.13 of \cite{HumCoxeter} is the basis for an argument in 
\S 4 of \cite{ErikssonDiscrete} showing that if $\lambda = w.\mu$ 
for $\mu \in -D$, then $\mu$ can be reached from $\lambda$ by a legal 
sequence.  Then we get the following 
characterization of the set of good positions: 

\noindent 
{\bf \TitsConeConvergenceResult\ (Eriksson)}\ \ 
{\sl The set of good positions for $(\Gamma,A)$ 
is precisely $-U$.} 

Our next result generalizes Remark 4.4 of \cite{Deodhar} to our 
current setting.  
This is needed for \HRTResult.  
For $J \subseteq I_{n}$, 
let $\Phi^{J} := \{\alpha \in \Phi^{+}\, |\, 
\alpha \not\in \mathrm{span}_{\mathbb{R}}\{\alpha_{j}\}_{j \in J}\}$. 

\noindent 
{\bf \DeodharProp}\ \ {\sl If $(\Gamma,A)$ is connected, 
$\Phi$ is infinite, and $J \subset I_{n}$ (proper), 
then $\Phi^{J}$ is infinite.}  

{\em Proof:} In the ``({\em ix}) $\Rightarrow$ ({\em ii})'' part of 
the proof of Proposition 4.2 in \cite{Deodhar}, 
assume $|\Phi^{J}| < \infty$ and begin reading at line -8 of page 
620.\hfill\QED

Proposition 3.2 of \cite{HRT} states 
that if $(\Gamma,A)$ is connected, $\sigma$ is standard, and $W$ 
is infinite, then $U_{0} = \{0\}$. 
In view of \HRTProofResults, we can use the proof of 
Proposition 3.2 of \cite{HRT} verbatim to get the 
generalization of that result stated as \HRTResult\ below.  
One can see that that 
proof will work if it is known that all 
$N(w)$ are finite; by \UnitalCorollary\ this is guaranteed by 
our hypothesis in the statement of \HRTResult\ requiring that 
$(\Gamma,A)$ is unital ON-cyclic.

\noindent 
{\bf \HRTResult}\ \ {\sl Suppose $(\Gamma,A)$ is connected and unital 
ON-cyclic and $W$ 
is infinite. Then 
$U_{0} = \{0\}$, i.e.\ $U$ is a ``strictly convex'' cone.} 
\hfill\QED
 
In contrast, for finite $W$ the overlap $U_{0} = U \cap (-U)$ is 
\underline{all} of $V^{*}$.  This is a consequence of the following 
result due to Vinberg (see \S 7 of \cite{Vinberg}).  The proof below uses 
numbers game reasoning. 

\noindent 
{\bf \TitsConeProp}\ \ 
{\sl If $W$ is finite, then $U = V^{*} = -U$, so $U_{0} = V^{*}$.} 

{\em Proof.} 
Since $W$ is finite, then by 
the Reduced Word 
Result it follows that the set of good positions is all of 
$V^{*}$.  
\TitsConeConvergenceResult\ now implies that $V^{*} = -U$, 
hence $U = V^{*}$ also.\hfill\QED 

When $(\Gamma,A)$ is connected and unital 
ON-cyclic, if a nonzero $\lambda \in D$ is good, then by 
\TheoremList,  
$W$ must be finite.  
This observation, together with the classification of finite Coxeter 
groups and reasoning based on the numbers game, 
is used in \S 6 of \cite{DonNumbers} to prove the following 
result, 
which we refer to in \VinbergRemark\ below as result ({\tt *}): 
{\sl If $(\Gamma,A)$ is 
connected, then $D \cap (-U) \not= \{0\}$ implies that 
$W$ is finite.}   
(We know of three proofs of statement ({\tt *}): 
See Theorem 6.1 of \cite{DonNumbers}; 
see \VinbergRemark\ below 
for a proof that uses \HRTResult, results borrowed from 
\cite{Vinberg}, and a classification result due to H.\ S.\ M.\ 
Coxeter; or see \S 4 of \cite{DE} for a proof 
that does not require \HRTResult\ or 
the classification of finite Coxeter groups.)  Now, it follows from 
the definitions that $U_{0} \not= \{0\}$ if and only if $D \cap (-U) 
\not= \{0\}$.  In view of \TitsConeProp, 
we thus obtain the following addition to the list of 
equivalences from Propositions 4.1 and 4.2 of \cite{Deodhar} for an 
irreducible Coxeter group to be finite: 

\noindent 
{\bf \TitsConeFiniteResult}\ \ {\sl Let $(\Gamma,A)$ be connected, 
so the Coxeter group $W$ is irreducible.  Then $W$ is finite if and 
only if $U_{0} \not= \{0\}$ if and only if $U_{0} = V^{*}$.}\hfill\QED  

See \S 2 of \cite{Krammer} for a proof of this result 
in the special case that 
the bilinear form $B$ for the representing space $V$ is 
symmetric. 

\noindent
{\bf \KumarRemark}\ \  
A Tits cone is similarly 
defined in the context of Kac--Moody theory 
e.g.\ \cite{Kac} \S 3.12, \cite{Kumar} \S 1.4.  
Let $A$ be a GCM.  
Here we follow Kac \cite{Kac} and the end of 
\S 2 above. 
The Kac--Moody Tits cone is the set $C := C_{A} := \{w.\lambda\, |\, 
w \in W, \lambda \in \mathfrak{h}_{\mathbb{R}} \mbox{ such that } 
\alpha_{i}(\lambda) \geq 0 \mbox{ for } 1 \leq i \leq n\} \subseteq 
\mathfrak{h}_{\mathbb{R}}$.  
When $A$ is nondegenerate ($\mbox{nullity}(A) = 0$), 
then $V = \mathfrak{h}_{\mathbb{R}}^{*}$ and hence $C$ and $U$ coincide.  
Now suppose $(\Gamma,A)$ is connected and $W$ is infinite. 
We have that the GCM graph $(\Gamma,A)$  
is unital ON-cyclic. 
Thus if  $A$ is nondegenerate,  
the result $C_{0} = \{0\}$ holds by \HRTResult.  
Allowing $\mbox{nullity}(A) \geq 0$, 
Kumar (personal communication) has supplied the 
following description of $C_{0}$: $C_{0} 
= \{v \in \mathfrak{h}_{\mathbb{R}}\, |\, \alpha_{i}(v) = 0 
\mbox{ for } 1 \leq i \leq n\}$.  
Here $\dim(C_{0}) = \mbox{nullity}(A)$.  
He notes that this statement may be deduced from 
Part (c) of Proposition 3.12 of \cite{Kac}.   
Note that the topological interior of the Kac--Moody Tits cone can 
never intersect its negative.  This follows from \cite{Kac} Exercise 
3.15 (see also \cite{Kumar} Exercise 1.4.E.1).\hfill\QED  

\noindent 
{\bf \VinbergRemark}\ \ 
In this remark we use \HRTResult\ and results from \cite{Vinberg} to 
prove the following version of result ({\tt *}) above: {\sl 
If $(\Gamma,A)$ is 
connected, then $W$ infinite implies that $U_{0} = \{0\}$.} (Then, 
Theorem 6.1 of \cite{DonNumbers} can be obtained as an easy consequence.)  
To prove this, we interpret our set-up here in terms of 
\cite{Vinberg}.   
Assume throughout this remark that $(\Gamma,A)$ is connected.  Our 
$V^{*}$ plays the role of 
Vinberg's $V$, our $D$ plays the role of his $K$, our $\alpha_{i}$'s play the 
role of his $\alpha_{i}$'s.  
Let us take elements $h_{i} \in 
V^{*}$ for which $h_{i}(\alpha_{j}) = a_{ij}$ for all $i, j \in 
I_{n}$.  Each $R_{i} := 
\sigma^{*}(s_{i}) \in GL(V^{*})$ is then a reflection in the 
sense of \S 1 of \cite{Vinberg}. From \cite{HumCoxeter} \S 5.13, we 
know that $w\! \stackrel{\circ}{D}  
\cap \stackrel{\circ}{D} = \emptyset$ for 
all $w \not= \varepsilon$ in $W$, where $\stackrel{\circ}{D}$ denotes 
the interior of $D$.  Therefore our $\sigma^{*}(W) \subset 
GL(V^{*})$ is, in Vinberg's language, a ``linear Coxeter group.'' 

With respect to some total ordering of $I_{n}$, 
think of an $n$-tuple $\velt = (v_{i})_{i \in I_{n}}$ as a column 
vector.  Denote by $\mathbf{0}$ the zero vector. 
For column vectors $\uelt, \velt$, say 
$\uelt > \velt$ (respectively $\uelt \geq \velt$) 
if $u_{i} > v_{i}$ (resp.\ $u_{i} \geq v_{i}$) for each $i \in I_{n}$.  
From \S 4 of \cite{Vinberg} or Chapter 4 of \cite{Kac}, 
we have that 
exactly one of the following three statements (+), (0), or (--) is 
true: (+) $\det A \not= 0$; there exists $\velt > \mathbf{0}$ such that 
$A\velt > \mathbf{0}$; $A\uelt \geq \mathbf{0}$ implies 
$\uelt > \mathbf{0}$ or $\uelt = \mathbf{0}$; 
(0) $\mbox{nullity}(A) = 1$;  there exists $\velt > \mathbf{0}$ such that 
$A\velt = \mathbf{0}$; $A\uelt = \mathbf{0}$ implies that $\uelt \geq 
\mathbf{0}$; (--) There exists $\velt > \mathbf{0}$ such that 
$A\velt < \mathbf{0}$; $A\uelt \geq \mathbf{0}$, $\uelt \geq 
\mathbf{0}$ imply that $\uelt = \mathbf{0}$.  Write $A=A^{+}$, $A= 
A^{0}$, or $A=A^{-}$ accordingly.  By Lemma 15 and Proposition 25 of 
\cite{Vinberg}, we see that $A=A^{+} \Rightarrow U=V^{*}$, $A=A^{0} 
\Rightarrow (\overline{U})_{0} = \mbox{span}(h_{i})_{i\in{I}_{n}}$, 
and $A=A^{-} \Rightarrow 
(\overline{U})_{0} = \{0\}$.  (Note that the 
set ``$\mbox{Ann}[\alpha]$'' in \cite{Vinberg} is $\{0\}$ here, 
since it is just $\{\lambda \in V^{*}\, |\, \lambda(v) = 0 \mbox{ for 
all } v \in V\}$.) 

Now we prove the version of result ({\tt *}) stated at the beginning of 
this remark.  We consider the three cases (+), (0), and (--).  First, 
suppose $A=A^{+}$.  Then by Proposition 22 of \cite{Vinberg}, $W$ 
must be finite, contrary to our hypothesis.  Second, suppose that 
$A=A^{0}$.  Then by Proposition 23 of \cite{Vinberg}, $W$ is an 
irreducible ``parabolic'' Coxeter group, also 
called an irreducible Euclidean reflection 
group, see e.g.\ \cite{Davis}.  The well-known classification of such 
groups is due to H.\ S.\ M.\ Coxeter \cite{Coxeter}. 
For our purposes, it is 
enough to observe that any such $(\Gamma,A)$ will possess a simple 
ON-cycle only in the case that $(\Gamma,A)$ itself is a simple 
ON-cycle with $m_{ij} = 3$ for any adjacent $\gamma_{i}$ and 
$\gamma_{j}$.  By Proposition 23 of \cite{Vinberg}, it follows that 
this ON-cycle is unital.  We conclude that whenever $A=A^{0}$, 
$(\Gamma,A)$ is unital ON-cyclic.  Since $W$ is infinite (by, say, 
Proposition 22 of \cite{Vinberg}), then by \HRTResult\ above, we 
have $U_{0} = \{0\}$.  Finally, if $A=A^{-}$, then $(\overline{U})_{0} = 
\{0\}$ implies that $U_{0} = \{0\}$. (Note that $W$ must infinite 
whenever $A=A^{-}$, 
by Proposition 22 of \cite{Vinberg}.)   
In any case, we see that when 
$(\Gamma,A)$ is connected and $W$ is infinite, then $U_{0} = 
\{0\}$.\hfill\QED   

\noindent 
{\bf \TitsConeExample}\ \ For $(\Gamma,A) =$ 
\setlength{\unitlength}{0.75in}
\begin{picture}(1.65,0.25)
\put(0.25,0){\begin{picture}(1,0)
            \put(0,0.1){\circle*{0.05}}
            \put(-0.20,-0.05){\large $\gamma_{1}$}
            \put(1,0.1){\circle*{0.05}}
            \put(1.05,-0.05){\large $\gamma_{2}$}
            \put(0,0.1){\line(1,0){1}}
            \put(0.2,0.1){\vector(1,0){0.1}}
            \put(0.8,0.1){\vector(-1,0){0.1}}
            \put(0.225,-0.05){\footnotesize $p$}
            \put(0.71,-0.05){\footnotesize $q$}
            \end{picture}}
\end{picture} with $pq = 4$, we have that $W$ is the 
infinite dihedral group.  Since $A=A^{0}$ in the notation of 
\VinbergRemark, then $(\overline{U})_{0} = 
\mbox{span}(h_{i})_{i\in\{1,2\}}$.  Relative to the basis 
$\{\omega_{i}\}_{i=1,2}$ for $V^{*}$ dual to the simple root basis 
$\{\alpha_{i}\}_{i=1,2}$ for $V$, we have $h_{1} = 2\omega_{1} - 
p\omega_{2}$ and $h_{2} = -q\omega_{1} + 2\omega_{2} = 
-\frac{q}{2}h_{1}$.  Then, $(\overline{U})_{0} = 
\mbox{span}(h_{i})_{i\in\{1,2\}} = 
\{x\omega_{1}+y\omega_{2}\, |\, y = -\frac{p}{2}x\}$.  In fact, using 
the computational approach of the proof of \TwoGeneratorAnalysis\ 
above, one can see that $U = \{x\omega_{1}+y\omega_{2}\, |\, 
y > -\frac{p}{2}x \mbox{ or } x=y=0\}$, and hence that $U_{0} = 
\{0\}$.\hfill\QED 

\noindent 
{\bf Acknowledgments}\ \ We thank Kimmo Eriksson for providing us with 
a copy of his thesis and for many helpful conversations during the 
preparation of this paper.  We thank Bob Proctor for his helpful 
feedback, which included the remarks concerning Weyl groups, and for 
sharing with us in advance some of the results of \cite{ProctorCoxeter}. 
We thank Shrawan Kumar for helpful comments concerning  Kac--Moody Tits 
cones.

\vspace*{-0.15in}
\renewcommand{\baselinestretch}{1}

\end{document}